\documentclass[a4paper,12pt,reqno,oneside]{amsart}

\usepackage{tikz}
\usetikzlibrary{decorations.pathmorphing}
\usetikzlibrary{cd}

\usepackage{setspace} 
\usepackage{typearea} 
\usepackage{amscd,amssymb,verbatim}
\usepackage
{graphicx} \setkeys{Gin}{width=\linewidth}
\usepackage{enumitem}
\usepackage{bbold}

\usepackage{xr-hyper}
\usepackage{cite} 
\usepackage[colorlinks,backref,final,bookmarksnumbered,bookmarks]{hyperref}
\usepackage[backrefs]{amsrefs}

\usepackage{ifthen}
\newcommand{\arxiv}[2][]{\ifthenelse{\equal{#1}{}}
{\href{http://arxiv.org/abs/#2}{\tt arXiv:#2}}
{\href{http://arxiv.org/abs/math/#2}{\tt arXiv:math.#1/#2}}}

\usepackage[T1,T2A]{fontenc}
\usepackage[utf8]{inputenc}

\makeatletter
\renewcommand\subsection{\@startsection
{subsection}{2}{0cm} 
{-\baselineskip}     
{0.5\baselineskip}   
{\sffamily}} 
\makeatother

\theoremstyle{plain}
\newtheorem{theorem}{Theorem}[section]

\newtheorem{corollary}[theorem]{Corollary}
\newtheorem{proposition}[theorem]{Proposition}
\newtheorem{problem}[theorem]{Problem}

\theoremstyle{definition}
\newtheorem{example}[theorem]{Example}

\newtheoremstyle{remark}
{}{}{}{}{\itshape}{}{ }{\thmname{#1}\thmnumber{ \itshape #2.}}
\theoremstyle{remark}
\newtheorem{remark}[theorem]{Remark}

\newtheoremstyle{concise}
{}{}{}{}{\bfseries}{}{ }{\thmnumber{#2.}\thmnote{ #3.}}
\theoremstyle{concise}

\DeclareFontEncoding{LS1}{}{}
\DeclareFontSubstitution{LS1}{stix}{m}{n}
\DeclareFontEncoding{LS2}{}{}
\DeclareFontSubstitution{LS2}{stix}{m}{n}

\DeclareMathVersion{language}
\DeclareMathVersion{metalanguage}

\def\formulas{\mathversion{language}}

\def\metameta{\mathversion{normal}}

\DeclareMathAlphabet{\mf}{OML}{zplm}{m}{it} 
\DeclareMathAlphabet{\ts}{OT1}{cmss}{m}{sl} 

\DeclareMathAlphabet{\fm}{U}{eur}{m}{n} 
\DeclareMathAlphabet{\tr}{OT1}{cmss}{m}{n} 
\DeclareMathAlphabet{\mm}{OML}{cmm}{m}{it} 

\DeclareMathAlphabet{\script}{LS1}{stixscr}{m}{n}
\DeclareMathAlphabet{\mathbbb}{U}{bbold}{m}{n}

\SetSymbolFont{letters}{language}{U}{eur}{m}{n}
\SetSymbolFont{letters}{metalanguage}{OML}{zplm}{m}{it}

\newcommand{\newsym}[5]{\fontfamily{#2}\fontencoding{#1}\fontseries{#3}\fontshape{#4}\selectfont\char#5}
\makeatletter
\newcommand{\newmathsymbol}[6]{#1{\@Pimathsymbol{#2}{#3}{#4}{#5}{#6}}}
\def\@Pimathsymbol#1#2#3#4#5{\mathchoice
  {\@Pim@thsymbol{#1}{#2}{#3}{#4}{#5}\tf@size}
  {\@Pim@thsymbol{#1}{#2}{#3}{#4}{#5}\tf@size}
  {\@Pim@thsymbol{#1}{#2}{#3}{#4}{#5}\sf@size}
  {\@Pim@thsymbol{#1}{#2}{#3}{#4}{#5}\ssf@size}}
\def\@Pim@thsymbol#1#2#3#4#5#6{\mbox{\fontsize{#6}{#6}\newsym{#1}{#2}{#3}{#4}{#5}}}
\makeatother

\def\too{\newmathsymbol{\mathrel}{LS1}{stixsf}{m}{n}{"99}}

\def\lgroup{\newmathsymbol{\mathopen}{LS2}{stixex}{m}{n}{"DC}}
\def\rgroup{\newmathsymbol{\mathclose}{LS2}{stixex}{m}{n}{"DD}}
\newcommand{\mq}[1]{\lgroup #1\rgroup\,}

\def\prin{\boldsymbol\cdot\hskip1.5pt}

\def\clbot{\bot}
\def\ab{{\mathchoice
  {\mbox{\larger[1]$\times$}}
  {\mbox{\larger[1]$\times$}}
  {\mbox{\larger[-2]$\times$}}
  {\mbox{\larger[-4]$\times$}}
}}
\def\triv{\checkmark}
\def\Bot{\newmathsymbol{\mathord}{U}{cmllr}{m}{n}{13}}
\def\Top{\raisebox{8pt}{\scalebox{1}[-1]{$\Bot$}}}

\def\KK{\hskip1pt\tikz[baseline=0pt,decoration={snake,amplitude=0.6pt,segment length=3pt}]
\draw[decorate] (0pt,-3pt) -- (0pt,9pt);\hskip1pt}

\def\from{\leftarrow}
\def\TO{\Rightarrow}
\def\To{\ \longrightarrow\ }
\def\tofrom{\leftrightarrow}
\def\Tofrom{\ \longleftrightarrow\ }
\def\imp{\Rightarrow}

\def\iff{\Leftrightarrow}

\def\turnstile{\vdash}
\def\Turnstile{\vDash}

\DeclareMathSymbol{\impord}{\mathord}{symbols}{41}
\DeclareMathSymbol{\mand}{\mathbin}{operators}{`\&}

\DeclareMathSymbol{\oc}{\mathord}{operators}{`!}
\DeclareMathSymbol{\wn}{\mathord}{operators}{`?}

\def\ocf{\newmathsymbol{\mathord}{OT1}{cmr}{m}{it}{"21}}
\def\wnf{\newmathsymbol{\mathord}{OT1}{cmr}{m}{it}{"3F}}

\DeclareMathSymbol{\alpha}     {\mathalpha}{letters}{"0B}
\DeclareMathSymbol{\beta}      {\mathalpha}{letters}{"0C}
\DeclareMathSymbol{\gamma}     {\mathalpha}{letters}{"0D}
\DeclareMathSymbol{\delta}     {\mathalpha}{letters}{"0E}
\DeclareMathSymbol{\epsilon}   {\mathalpha}{letters}{"0F}
\DeclareMathSymbol{\zeta}      {\mathalpha}{letters}{"10}
\DeclareMathSymbol{\eta}       {\mathalpha}{letters}{"11}
\DeclareMathSymbol{\theta}     {\mathalpha}{letters}{"12}
\DeclareMathSymbol{\iota}      {\mathalpha}{letters}{"13}
\DeclareMathSymbol{\kappa}     {\mathalpha}{letters}{"14}
\DeclareMathSymbol{\lambda}    {\mathalpha}{letters}{"15}
\DeclareMathSymbol{\mu}        {\mathalpha}{letters}{"16}
\DeclareMathSymbol{\nu}        {\mathalpha}{letters}{"17}
\DeclareMathSymbol{\xi}        {\mathalpha}{letters}{"18}
\DeclareMathSymbol{\pi}        {\mathalpha}{letters}{"19}
\DeclareMathSymbol{\rho}       {\mathalpha}{letters}{"1A}
\DeclareMathSymbol{\sigma}     {\mathalpha}{letters}{"1B}
\DeclareMathSymbol{\tau}       {\mathalpha}{letters}{"1C}
\DeclareMathSymbol{\upsilon}   {\mathalpha}{letters}{"1D}
\DeclareMathSymbol{\phi}       {\mathalpha}{letters}{"1E}
\DeclareMathSymbol{\chi}       {\mathalpha}{letters}{"1F}
\DeclareMathSymbol{\psi}       {\mathalpha}{letters}{"20}
\DeclareMathSymbol{\omega}     {\mathalpha}{letters}{"21}
\DeclareMathSymbol{\varepsilon}{\mathalpha}{letters}{"22}
\DeclareMathSymbol{\vartheta}  {\mathalpha}{letters}{"23}
\DeclareMathSymbol{\varpi}     {\mathalpha}{letters}{"24}
\DeclareMathSymbol{\varrho}    {\mathalpha}{letters}{"25}
\DeclareMathSymbol{\varsigma}  {\mathalpha}{letters}{"26}
\DeclareMathSymbol{\varphi}    {\mathalpha}{letters}{"27}

\DeclareMathSymbol{\bbomega}{\mathalpha}{letters}{"7F}

\def\C{\mathbbb{C}}
\def\D{\script{D}}
\def\Ds{\mathcal{D}}
\def\Dec{\mathfrak{D}}
\def\E{\script{E}}
\def\F{\mathcal{F}}
\def\Fm{\script{F}}
\def\G{\mathcal{G}}
\def\H{\mathcal{H}}

\def\N{\mathbbb{N}}
\def\O{\script{O}}
\def\Q{\mathbbb{Q}}
\def\Qm{\script{Q}}
\def\R{\mathbbb{R}}
\def\Rt{\mathbbb{\omega}}

\def\Stab{\mathfrak{S}}
\def\T{\script{T}}

\def\U{\mathcal{U}}
\def\V{\mathcal{V}}

\def\0{\mathbbb{0}}
\def\1{\mathbbb{1}}
\def\m{\mathbbb{\mu}}
\def\q{\mathbbb{q}}
\def\iass{\script{u}}
\def\pval{\script{v}}

\def\x{\times}
\def\but{\setminus}
\def\emb{\hookrightarrow}

\def\phi{\varphi}
\def\emptyset{\varnothing}
\def\xr#1{\xrightarrow{#1}}

\renewcommand{\:}{\colon}

\DeclareMathOperator{\id}{id}
\DeclareMathOperator{\Int}{Int}
\DeclareMathOperator{\Cl}{Cl}

\DeclareMathOperator{\Hom}{Hom}
\DeclareMathOperator{\Supp}{Supp}
\DeclareMathOperator{\Char}{Char}
\DeclareMathOperator{\Dom}{Dom}

\DeclareMathOperator{\Bew}{Bew}

\def\octo{\oc_{_\to}} \def\wnto{\wn_{_\to}}

\def\o#1{\#\oldstylenums{#1}}

\begin{document}

\title[A Galois connection betw. intuitionistic \& classical logics. II: Semantics]
{A Galois connection between intuitionistic and classical logics. II: Semantics}
\author{Sergey A. Melikhov}
\address{Steklov Mathematical Institute of Russian Academy of Sciences,
ul.\ Gubkina 8, Moscow, 119991 Russia}
\email{melikhov@mi-ras.ru}
\thanks{Supported by The Fund for Math and
Russian Foundation for Basic Research Grant No.\ 15-01-06302}

\begin{abstract}
Three classes of models of QHC, the predicate joint logic of problems and propositions, are constructed, 
including a class of subset/sheaf-valued models that is related to solutions of some actual problems (such as 
solutions of algebraic equations).

To test the models, we consider a number of principles and rules, which empirically appear to cover all 
``sufficiently simple'' natural conjectures about the behaviour of the operators $\oc$ and $\wn$, and 
include two hypotheses put forward by Hilbert and Kolmogorov, as formalized in the language of QHC.
Each of these turns out to be either derivable in QHC or equivalent to one of 14 principles and rules, 
of which 11 are conservative over classical and intuitionistic logics.
The three classes of models together suffice to confirm the independence of these 11 principles and rules, 
and to determine the full lattice of implications between them, apart from one potential implication.
\end{abstract}

\maketitle

\section{Introduction}

The paper is organized as follows.
Models of QHC are discussed in the last chapter (\S\ref{models}), which is rather short due to
the previous preparation.
Namely, since equivalent principles and rules have the same interpretation in every model,
our study of models is aided by a preceding study of the lattice of implications between principles 
and rules (also known as the lattice of extensions of QHC) in \S\ref{principles section}.
Some of these principles and rules are also of independent interest from the viewpoint of informal 
semantics of QHC.
We also look at both formal and informal semantics of the meta-logic of QHC 
in \S\ref{semantics of meta-logic}.

In the following section (\S\ref{equations}) we discuss how the subset/sheaf-valued models of
\S\ref{sheaf model} arise naturally in elementary mathematics.

\subsection{Algebraic equations}\label{equations}

Let $f\:X\to B$ be a continuous map (for instance, this could be a real polynomial $f\:\R\to\R$ or a complex
polynomial $f\:\C\to\C$, $f(x)=a_nx^n+\dots+a_0$), and let us consider the following parametric problem
$\Gamma_f(b)$, $b\in B$:
\begin{center}
{\it Find a solution of the equation $f(x)=b$.}
\end{center}

If $g\:Y\to B$ is another continuous map, then $\Gamma_f(b)\land\Gamma_g(b)$ is the problem
of finding a solution of the system of equations
$$\begin{cases}
f(x)=b\\
g(y)=b.
\end{cases}$$
Similarly, $\Gamma_f(b)\lor\Gamma_g(b)$ is the problem of finding a solution of the union of equations%
\footnote{The square bracket notation for union of equations (and union of systems of equations) is commonly
used in some cultures, e.g.\ in Russia and Ukraine (for instance, when dealing with absolute value equations, 
equations reducible to the form $f(x)g(x)=0$, etc.).}
$$\left[\begin{gathered}
f(x)=b\\
g(y)=b.
\end{gathered}
\right.$$
(Thus the traditional signs $\{$ and $\left[\right.$ can thought of as notationally representing
the intuitionistic, and not classical, conjunction and disjunction, as long as we understand
the line containing an equation as implicitly asking one to find a solution of this equation.)

The proposition $\wn\Gamma_f(b)$ reads:
\begin{center}
{\it The equation $f(x)=b$ has a solution,}
\end{center}
This is a parametric proposition with parameter $b$ running over $B$; in other words, a unary predicate 
on $B$.
Its {\it truth value} $|\wn\Gamma_f|$ is the set of all $b\in B$ for which $\wn\Gamma_f(b)$ is true.
This set is nothing but $f(X)$.

The truth value of the conjunction $\wn\Gamma_f\land\wn\Gamma_g$ is
the intersection $f(X)\cap g(Y)$;
and the truth value of the disjunction $\wn\Gamma_f\lor\wn\Gamma_g$ is the union $f(X)\cup g(X)$.

Each point-inverse $f^{-1}(b)$ can be identified with the set of solutions of $\Gamma_f(b)$
(this identification amounts to the hypothesis of ``functional extensionality'').

The solutions of $\Gamma_f(b)\lor\Gamma_g(b)$ are the elements of $f^{-1}(b)\sqcup g^{-1}(b)$,
which are in turn identified with the solutions of $\Gamma_{f\sqcup g}(b)$, where $f\sqcup g$
denotes the obvious map to $B$ from the disjoint union $X\sqcup Y$.
The solutions of $\Gamma_f(b)\land\Gamma_g(b)$ are the elements of
$f^{-1}(b)\x g^{-1}(b)$, which are in turn identified with the solutions of $\Gamma_{f\x g}(b)$,
where $f\x g$ denotes the obvious map to $B$ from the fiberwise product
$X\x_B Y=\{(x,y)\in X\x Y\mid f(x)=g(y)\}$.

Thus we can make the following identifications of the problems with parameter $b$ running over $B$:
\[\Gamma_f\land\Gamma_g=\Gamma_{f\x g}\quad\text{and}\quad
\Gamma_f\lor\Gamma_g=\Gamma_{f\sqcup g}.\]
Let us note that if $f$ and $g$ happen to be sheaves of sets, then $f\x g$ and $f\sqcup g$ are their usual
product and coproduct.

We see that, in our setting of parametric equations, $\wn$ has the effect of converting categorical notions
into set-theoretic ones: the image of the product map, $f\x g$, is the intersection $f(X)\cap g(Y)$; and
the image of the coproduct map, $f\sqcup g$, is the union $f(X)\cup g(X)$.
In more syntactic terms,
\[|\wn\Gamma_{f\x g}|=|\wn(\Gamma_f\land\Gamma_g)|=
|\wn\Gamma_f\land\wn\Gamma_g|=|\wn\Gamma_f|\cap|\wn\Gamma_g|;\]
\[|\wn\Gamma_{f\sqcup g}|=|\wn(\Gamma_f\lor\Gamma_g)|=
|\wn\Gamma_f\lor\wn\Gamma_g|=|\wn\Gamma_f|\cup|\wn\Gamma_g|.\]

Now let us consider the parametric problem $\Gamma_g(b)\to\Gamma_f(b)$ of reducing the equation
$f(x)=b$ to the equation $g(y)=b$.
What is a solution of this problem, and when does one exist?
One thing is immediately clear: if there is a reduction of $f(x)=b$ to $g(y)=b$, and the latter equation has
a solution, the former must also have one; in symbols,
$\wn\big(\Gamma_g(b)\to\Gamma_f(b)\big)\to\big(\wn\Gamma_g(b)\to\wn\Gamma_f(b)\big)$.
Let us note, however, that the implication ``if $f(x)=b$ has a solution, then so does $g(x)=b$'' is
a vacuous condition if $f$ and $g$ are complex polynomials of nonzero degree (since both equations
necessarily have solutions by the fundamental theorem of algebra).
But this does not mean that it is meaningless to ask which polynomial equations reduce to which.

Indeed, everybody knows how to solve the quadratic equation $x^2+px=b$.
Upon rewriting it in the form $(x+\frac{p}{2})^2-\frac{p^2}4=b$ we see that it {\it reduces} to
the equation $y^2-\frac{p^2}4=b$, which is suddenly easy to solve.
Here $y=x+\frac p2$, but the actual {\it reduction} of the former equation to the latter one is given by
the substitution $x=y-\frac{p}2$, whose result $(y-\frac{p}2)^2+p(y-\frac{p}2)=b$ ``miraculously''
lacks the linear term.
In general, an equation $f(x)=b$ reduces trivially to any equation of the form $f(\phi(y))=b$: by finding
a solution $y_0$ of the latter equation, we immediately obtain the solution $x_0=\phi(y_0)$ of
the former one.
To find a substantial reduction, one seeks to rewrite the given equation $f(x)=b$ in the form
$g(\psi(x))=b$ so that the equation $\psi(x)=y$ can be solved for $x$; this solution yields
the desired function $\phi=\psi^{-1}$.%
\footnote{\label{Vieta}
An illustration of the power of this method is Vieta's solution of the cubic equation $x^3+ax^2+cx=b$.
Firstly, the substitution $x=y-\frac{a}3$ reduces it to the form $y^3+py+q=b$.
Secondly, the substitution $y=z-\frac{p}{3z}$ results in the equation
$z^3+q-\frac{p^3}{27z^3}=b$, as long as $z\ne 0$.
Thirdly, the substitution $z=\sqrt[3]{w}$ yields a reduction to the equation $w+q-\frac{p^3}{27w}=b$,
which upon multiplication by $w$ becomes quadratic (and so can be solved as above).}

Abstracting from the algebraic setting, let us say that a {\it reduction} of the parametric equation
$f(x)=b$ to the parametric equation $g(y)=b$ assigns to every solution $y_0$ of the latter one
a solution $x_0$ of the former one, by an assignment $\phi$ that depends continuously on $b$.
(In the example above, $\phi$ does not depend on $b$ at all.)
Thus it is a family of maps $g^{-1}(b)\to f^{-1}(b)$ depending continuously on $b$; or, in other words,
a continuous map $\phi$ making the following diagram commute:
\[\begin{tikzcd}[ampersand replacement=\&,column sep=large]
Y\ar{rd}[swap]{g} \ar{r}{\phi}
\& X\ar{d}{f}\\
\& B
\end{tikzcd}\]

However, if we return to our quadratic equation $x^2+px=b$, we have only seen that it reduces to
$y^2-\frac{p^2}4=b$.
To solve the latter, we would first reduce it to the linear equation $z-\frac{p^2}4=b$ via the substitution
$y=\pm\sqrt z$, and then the linear equation can be solved directly.
Here the reducing ``function'' $\phi=\psi^{-1}$ is multivalued and, if we work over $\R$, only partially
defined.%
\footnote{Vieta's last substitution, $z=\sqrt[3]{w}$, is also multivalued if we work over $\C$.}
Thus if $b>\frac{p^2}4$, the quadratic equation has no real solutions, so in this case we do not really have
{\it any} reduction of the quadratic equation to the linear one (which has one solution for all $b$).
If $b<\frac{p^2}4$, the quadratic equation has two solutions, so we actually have two {\it distinct}
reductions in this case.

Thus algebraic substitutions actually yield reductions only locally.
Hence, let us say that a {\it local reduction} of $f(x)=*$ to $g(y)=*$ is a continuous map from
the preimage $g^{-1}(U)$ of some open set $U\subset B$ into $f^{-1}(U)$ making the following diagram commute:
\[\begin{tikzcd}[ampersand replacement=\&,column sep=large]
g^{-1}(U)\ar{rd}[swap]{g} \ar{r}{\phi}
\& f^{-1}(U)\ar{d}{f}\\
\& U.
\end{tikzcd}\]
For example, our quadratic equation has two local reductions to the linear equation over $\R$,
both defined over $U=(-\infty,\frac{p^2}4)$; they are given by $y=\sqrt z$ and $y=-\sqrt z$.

Now just like the parametric problems $\Gamma_f\land\Gamma_g$ and $\Gamma_f\sqcup\Gamma_g$,
the implication between the problem forms $\Gamma_g\to\Gamma_f$ must be a parametric problem.
We can now define a solution of $(\Gamma_g\to\Gamma_f)(b)$ to be the {\it germ} $\phi_b$ at $b$ of
a local reduction $\phi$ of $f(x)=*$ to $g(y)=*$, defined over some neighborhood $U$ of $b$.

If $\Hom(g,f)_b$ denotes the set of all such germs $\phi_b$ at a fixed point $b\in B$, there is
a natural topology on $\bigcup_{b\in B}\Hom(g,f)_b$, where a neighborhood of $\phi_b$ is given by
the germs $\phi_{b'}$ for all $b'\in U$, where $U$ is a neighborhood of $b$.
The map $\Hom(g,f)\:\bigcup_{b\in B}\Hom(g,f)_b\to B$, $\phi_b\mapsto b$, is known as
the {\it sheaf of germs of maps $Y\to X$ over $B$}, and we conclude that the problem form
$\Gamma_g\to\Gamma_f$ can be identified with $\Gamma_{\Hom(g,f)}$.

Now there seems to be a little issue with this definition.
There is the trivial equation $y=b$, corresponding to the case $Y=B$, $g=\id\:B\to B$.
If the problem $\Gamma_{\id}$ of finding a solution of this trivial equation is really a trivial problem,
then any problem $\Gamma_f$ should be equivalent to the problem $\Gamma_{\id}\to\Gamma_f$
of reducing the equation $f(x)=b$ to the trivial equation $y=b$.
But $\Hom(\id,f)$ is also known as the {\it sheaf of germs of sections} of the map $f$, and $\Gamma_{\Hom(\id,f)}$
is inequivalent to $\Gamma_f$ already for polynomial $f$, either complex (due to multiple roots) or real
(due to roots of even multiplicity); see \cite{M0}*{Example \ref{int:equations}} for the details.

Our approach to this issue is that it is not the definition of reduction between parametric problems
that is not quite right, but the very formulation of our problem $\Gamma_f$.
In fact, if we replace it by the problem $\Delta_f:=\Gamma_{\Hom(\id,f)}$, the issue disappears.
In closed terms, the problem $\Delta_f(b)$ asks to find only {\it stable solutions} of
the equation $f(x)=b$, that is, such $x_0\in X$ that $f(x_0)=b$ and there exists a neighborhood $U$
of $b$ in $B$ over which $f$ has a section (that is, a map $s\:U\to X$ such that the composition
$U\xr{s}X\xr{f}B$ equals the inclusion map $U\emb B$).

What has been said above about parametric problems of the form $\Gamma_f$ carries over to those of
the form $\Delta_f$; the result is nothing but a special case of the subset/sheaf-valued models 
of \S\ref{sheaf model}.

\section{Semantics of the meta-logic} \label{semantics of meta-logic}

What follows is a review of \cite{M0} as applied to the setup of QHC, and also some minor elaborations on it.

\subsection{Tarski-style semantics}

The types $\0$, $\1_c$, $\1_i$ and $\m$ are interpreted by sets $\D$, $\O_c$, $\O_i$ and $\Qm$. 
A map $\wnf\:\Qm\to\{\Top,\Bot\}$ is fixed.

The reflection operators $\Rt_c\:\1_c\too\m$ and $\Rt_i\:\1_i\too\m$ are interpreted by maps
$\ocf_c\:\O_c\to\Qm$ and $\ocf_i\:\O_i\to\Qm$.
The meta-connectives are interpreted by maps $|{\mand}|,|{\imp}|\:\Qm\x\Qm\to\Qm$; the first-order meta-quantifier
by a map $|\q|\:\Hom(\D,\Qm)\to\Qm$; and the second-order meta-quantifiers by  
$|\q_i^n|\:\Hom(\Hom(\D^n,\O_i),\Qm)\to\Qm$ and $|\q_c^n|\:\Hom(\Hom(\D^n,\O_c),\Qm)\to\Qm$.

The conversion operators are interpreted by maps $|\oc|\:\O_c\to\O_i$ and $|\wn|\:\O_i\to\O_c$.
Binary connectives are interpreted by maps $|{\land}|,|{\lor}|,|{\to}|\:\O_i\x\O_i\to\O_i$ (intuitionistic) and 
$|{\land}|,|{\lor}|,|{\to}|\:\O_c\x\O_c\to\O_c$ (classical); unary by maps $|\neg|\:\O_i\to\O_i$ (intuitionistic) 
and $|\neg|\:\O_c\to\O_c$ (classical); and nullary by elements $|\top|,|\bot|\in\O_i$ and 
$|\triv|,|\ab|\in\O_c$.
Intuitionistic quantifiers are interpreted by maps $|\forall|,|\exists|\:\Hom(\D,\O_i)\to\O_i$ and classical 
by maps $|\forall|,|\exists|\:\Hom(\D,\O_c)\to\O_c$, where $\Hom(X,Y)$ denotes the set of maps $X\to Y$.

The above is called a {\it meta-interpretation} of the language of QHC.
When $\Qm=\{\Top,\Bot\}$, $\wnf$ is the identity map, and the meta-connectives and meta-quantifiers are
interpreted according to the usual truth tables, the meta-interpretation is called {\it two-valued}, or
simply an {\it interpretation}.
Interpretations are, in a certain sense, models of the meta-logic (for a more precise formulation see 
\cite{M0}*{\S\ref{int:models}}).
In the present paper, we use only interpretations to study formal semantics of QHC.
However, for motivational purposes we will also employ a certain ``natural'' contentual semantics of QHC
(described in \S\ref{natural} below), which is based on a many-valued meta-interpretation.
It should also be noted that some realizability-type ``models'' of intuitionistic logic are in fact
meta-interpretations and not interpretations 
(see \cite{M0}*{\S\ref{int:uniform realizability}}).

We recall that an {\it $n$-formula} is a $\lambda$-expression of the form $x_1,\dots,x_n\mapsto F$, where $F$ 
is a formula and $x_1,\dots,x_n$ are pairwise distinct individual variables.
A $\lambda$-expression of the form $p_1,\dots,p_l,\gamma_1,\dots,\gamma_m\mapsto G$, where $G$ is an 
$n$-formula, $p_1,\dots,p_l$ are pairwise distinct predicate variables of (not necessarily distinct) arities
$q_1,\dots,q_l$ and $\gamma_1,\dots,\gamma_m$ are pairwise distinct problem variables of not (necessarily 
distinct) arities $r_1,\dots,r_m$, is called an {\it $(n,\vec q,\vec r)$-formula}, where
$\vec q=(q_1,\dots,q_l)$ and $\vec r=(r_1,\dots,r_m)$.
It can also be called an {\it $(n,\vec q,\vec r)$-c-formula} or an {\it $(n,\vec q,\vec r)$-i-formula} if $G$ 
is an $n$-c-formula or an $n$-i-formula.
One can similarly define $n$-meta-formulas, which are of type $\0^n\to\m$, and $(n,\vec q,\vec r)$-meta-formulas,
which are of type 
\[\0^n\x(\0^{q_1}\too\1)\x\dots\x(\0^{q_l}\too\1)\x(\0^{r_1}\too\1)\x\dots\x(\0^{r_m}\too\1)\too\m,\] where 
$\vec q=(q_1,\dots,q_l)$ and $\vec r=(r_1,\dots,r_m)$.

Given a meta-interpretation, a $\lambda$-closed $(n,\vec q,\vec r)$-c-formula $F$ is interpreted by 
a straightforward recursion by a map $|F|\:\D^n\x\Hom(\D^{\vec q,\vec r},\O)\to\O_c$, where 
$\Hom(\D^{\vec q,\vec r},\O)$ denotes
\[\Hom(\D^{q_1},\O_c)\x\dots\x\Hom(\D^{q_l},\O_c)\x\Hom(\D^{r_1},\O_i)\x\dots\x\Hom(\D^{r_m},\O_i),\]
where $\vec q=(q_1,\dots,q_l)$ and $\vec r=(r_1,\dots,r_m)$.
Similarly, a $\lambda$-closed $(n,\vec q,\vec r)$-i-formula $\Phi$ is interpreted by a map 
$|\Phi|\:\D^n\x\Hom(\D^{\vec q,\vec r},\O)\to\O_i$.
Using this, a $\lambda$-closed $(n,\vec q,\vec r)$-meta-formula $\F$ is interpreted by a map 
$|\F|\:\D^n\x\Hom(\D^{\vec q,\vec r},\O_i)\to\Qm$.

In particular, every $\lambda$-closed meta-formula $\F$ is interpreted by an element $|\F|\in\Qm$.
A $\lambda$-closed meta-formula $\F$ is called {\it valid} in a given meta-interpretation if $\wnf(|\F|)=\Top$.
If $\Ds$ is a derivation system of QHC, an interpretation $M$ of the language of QHC is called a {\it model}
of QHC if the $\lambda$-closed meta-formula $\Ds$ is valid in $M$. 

\subsection{Valuations and valuation fields}

A {\it valuation} $\pval$ assigns an element $\pval(p)\in\Hom(\D^n,\O_c)$ to each $n$-ary
predicate variable $p$, and an element $\pval(\gamma)\in\Hom(\D^n,\O_i)$ to each $n$-ary problem variable
$\gamma$, for each $n=0,1,2,\dots$.
Given a valuation $\pval$, every closed c-formula $F$ is interpreted by an element 
$|F|^\pval\in\O_c$, every closed i-formula $\Phi$ by an element $|\Phi|^\pval\in\O_i$, and every closed 
meta-formula $\F$ by an element $|\F|^\pval\in\Qm$.

An {\it assignment} $\iass$ assigns an element $\iass(x)\in\D$ to each individual variable $x$.
Given a valuation $\pval$ and an assignment $\iass$, every c-formula $F$ is interpreted 
by an element $|F|^\pval_\iass\in\O_c$, every i-formula $\Phi$ by an element $|\Phi|^\pval_\iass\in\O_i$, 
and every meta-formula $\F$ by an element $|\F|^\pval_\iass\in\Qm$.

Given a meta-interpretation of the language of QHC, let us consider a family $\H$ of subsets 
$\H_n^c\subset\Hom(\D^n,\O_c)$ and $\H_n^i\subset\Hom(\D^n,\O_i)$, $n=0,1,2,\dots$.
The family $\H$ is called a {\it valuation field} if it is closed with respect to the interpretations of 
connectives, quantifiers and the conversion operators.
A valuation $\pval$ is said to be {\it contained} in $\H$ if $\pval(p)\in\H_n^c$ and $\pval(\phi)\in\H_n^i$ 
for every $n$-ary predicate variable $p$ and every $n$-ary problem variable $\phi$, $n=0,1,2,\dots$.
The valuation field {\it generated by} a valuation $\pval$, denoted $\left<\pval\right>$, is the smallest 
valuation field that contains $\pval$.
Clearly, $\left<\pval\right>=(\H_n^c,\H_n^i)$, where each $\H_n^c$ consists of all $|F|^\pval$ such that 
$F$ is an $n$-c-formula, and each $\H_n^i$ consists of all $|\Phi|^\pval$ such that $\Phi$ is an $n$-i-formula.

Given an interpretation of the language of QHC, its {\it restriction} over a valuation field $\H$ is obtained
by re-interpreting the second-order meta-quantifiers by functions $|\q^n_c|\:\Hom(\H_n^c,\Qm\big)\to\Qm$
and $|\q^n_i|\:\Hom(\H_n^i,\Qm\big)\to\Qm$ defined in the same way as before, i.e.\ by
$|\q^n_c|(f)=\Top$ if and only $f(p)=\Top$ for all $p\in\H_n^c$ and by $|\q^n_i|(f)=\Top$ if and only 
$f(\phi)=\Top$ for all $\phi\in\H_n^i$.
Interpretations restricted over valuation fields are, in a certain sense, models of the meta-logic 
(for a more precise formulation see \cite{M0}*{\S\ref{int:models} and \S\ref{int:valuation fields}}).
An interpretation restricted over a valuation field $\H$ is called a {\it model} of QHC {\it with respect 
to $\H$} if some derivation system $\Ds$ of QHC is valid in the restricted interpretation.

\subsection{Clarified BHK interpretation} \label{clarified}

Kolmogorov's problem interpretation of intuitionistic logic \cite{Kol} explains not only logical, but, in 
disguise, also meta-logical connectives and quantifiers in the same language of solutions of problems.
The logical part was discussed in \cite{M1}*{\S\ref{g1:deductive1}}, and the meta-logical part will be discussed 
now.

In fact, we will follow a slight modification of Kolmogorov's original interpretation, called 
the ``clarified'' or (in extended form discussed below) ``meta-clarified'' BHK interpretation in \cite{M0}.
Namely, Kolmogorov's interpretation of, say, $\turnstile\fm\gamma\lor\neg\fm\gamma$ is that it is 
``the problem of finding a general method of solving the problem $\Gamma\lor\neg\Gamma$ for every 
contentful (e.g.\ mathematical) problem $\Gamma$''.
But in our approach, the informal meta-meta-logic of judgements about formal derivations in the meta-logic
is classical; and it does not seem to be appropriate to interpret such judgements by problems. 

Thus according to the clarified BHK interpretation,
\begin{enumerate} 
\item the principle $\prin\fm\gamma\lor\neg\fm\gamma$ is interpreted by the very same (meta-)problem mentioned 
above, {\it the problem of finding a general method of solving the problem $\Gamma\lor\neg\Gamma$ for every 
contentful problem $\Gamma$}; 
\item the judgement $\turnstile\prin\fm\gamma\lor\neg\fm\gamma$ asserting that the latter principle is
derivable is interpreted by the judgement $\Turnstile\prin\fm\gamma\lor\neg\fm\gamma$ that the said 
meta-problem has a solution.
\end{enumerate}

Let us note that either version includes a constructive quantification (``general method'') over all
concrete problems.
This does not seem to be compatible with usual, Tarski-style model theory, which only deals with 
classical quantification (e.g.\ over all valuations).
However, it is compatible with meta-interpretations as defined above, including some specific ones, such as
absolute realizability \cite{M0}*{\S\ref{int:uniform realizability}}.

\subsection{Kolmogorov's semantic consequence}

All interpretations considered above are either traditional, ``Tarski-style'' interpretations, which are based,
in particular, on Tarski's concept of semantic consequence; or meta-interpretations, which can be called 
``generalized Tarski-style''.

We will next do something radically different, and consider alternative, ``Kolmogorov-style'' interpretations
and meta-interpretations, which are based, in particular, on Kolmogorov's concept of semantic consequence, 
which is also found in his problem interpretation of intuitionistic logic \cite{Kol}.

Kolmogorov's problem interpretation uses the same language of solutions of problems to explain both logical
and (in disguise) meta-logical connectives and quantifiers.
It is this coincidence of $\O$ and $\Qm$ that makes Kolmogorov's concept of semantic consequence particularly 
apposite to his context.

Here is Kolmogorov's interpretation of the {\it modus ponens} rule, 
\[\fm{\frac{\phi,\phi\to\psi}{\psi}},\]
in the zero-order case which he considered.
We will break it into steps, which will be presented using our meta-logical notation; but upon eliminating 
the notation it is essentially Kolmogorov's original interpretation (see \cite{M0}*{\S\ref{int:K-consequence}}
for a more detailed discussion). 
\begin{enumerate}
\item The above rule is interpreted as the problem of finding a general method of solving a certain problem 
\[|{\prin\Phi\ \mand\ \prin\Phi\to\Psi\ \imp\ \prin\Psi}|\] for any given formulas 
$\Phi$ and $\Psi$.
\item A solution of the latter problem is in turn a general method of converting any given solutions of
$|{\prin\Phi}|$ and $|{\prin\Phi\to\Psi}|$ into a solution of $|{\prin\Psi}|$.
\item In the zero-order case, $\prin\Phi$, $\prin\Psi$ and $\prin\Phi\to\Psi$ are abbreviations for 
$\mq{\vec\gamma}\Phi$,  $\mq{\vec\gamma}\Psi$ and $\mq{\vec\gamma}\Phi\to\Psi$, where 
$\vec\gamma=(\gamma_1,\dots,\gamma_m)$ is the tuple of all problem variables that occur in $\Phi$ and $\Psi$.
A solution of e.g.\ $|\mq{\gamma}\Phi|$ is a general method of solving the contentful problem 
$|{\vec\gamma\mapsto\Phi}|(\vec\Gamma)$ for any tuple $\vec\Gamma=(\Gamma_1,\dots,\Gamma_m)$ of 
contentful (e.g.\ mathematical) problems.
\item What is a solution of $|{\vec\gamma\mapsto\Phi}|(\vec\Gamma)$, or 
$|{\vec\gamma\mapsto\Psi}|(\vec\Gamma)$, or 
$|{\vec\gamma\mapsto\Phi\to\Psi}|(\vec\Gamma)$ is determined by the usual BHK interpretation 
(see \cite{M1}*{\S\ref{g1:deductive1}}).
\end{enumerate}

Kolmogorov offered similar interpretations also for two other rules, and clearly this approach applies to 
an arbitrary rule just as well, including those in first-order logic (see details in 
\cite{M0}*{\S\ref{int:K-consequence}}).
Thus we will speak of {\it Kolmogorov's interpretation of rules}.

An obvious issue with Kolmogorov's interpretation of rules is that it applies what in essence is constructive 
quantification (``general method'') to syntactic entities (``for any given formulas''), thereby creating 
a dangerous mix of syntax and semantics.

But it turns out that this mix of syntax and semantics can be disentangled.
What is important in Kolmogorov's interpretation of rules is that the formulas being quantified over
(such as $\Phi$ and $\Psi$ in the above example) are treated as functions of the problem variables
$\vec\gamma$ (and in the first order case also of the individual variables $\vec x$).
But such a functional dependence can well be modelled on the semantic level.

\subsection{Alternative semantics}

Here we adopt the alternative interpretation of the meta-logic \cite{M0}*{\S\ref{int:Frege-style}}
(see also further details in \cite{M0}*{\S\ref{int:Frege-style2}}) to the setup of QHC.
The notation is somewhat different from that in \cite{M0}*{\S\ref{int:Frege-style}} 
(more natural, but less concise).

If $N$ is a finite set of individual variables, a {\it semantic $N$-term} is a function $T\:\D^{\#N}\to\D$,
where $\#N$ denotes the cardinality of $N$. 
A {\it semantic term} is an element $(N,T)$ of the disjoint union 
$|\T|:=\bigsqcup_{N\in\U_*}\Hom(\D^{\#N},\D)$, where $\U$ is the set of all individual variables and
$\U_*$ is the set of its finite subsets.

If additionally $M^c$ is a finite set of predicate variables and $M^i$ is a finite set of problem variables, 
and x stands for either c or i, a {\it semantic $(N,M^c,M^i,k)$-x-formula} is a function 
$\Phi\:\D^{\#N}\x\Hom(\D^{\#M^c,\#M^i},\O)\to\Hom(\D^k,\O_x)$,
where $\#M^x$ denotes the vector of arities of the elements of $M^x$.
A {\it semantic $k$-x-formula} is an element $(N,M^c,M^i,\Phi)$ of the disjoint union
$|\Fm_k^x|:=\bigsqcup_{(N,M^c,M^i)\in\U_*\x\V^c_*\x\V^i_*}\Hom\big(\D^{\#N}\x\Hom(\D^{\#M^c,\#M_i},\O),
\Hom(\D^k,\O_x)\big)$, where $\V^c$ is the set of all predicate variables,
$\V^i$ is the set of all problem variables, and $\V^x_*$ is the set of all finite subsets of $\V^x$.
We will also denote $|\Fm_{q_1}^c|\x\dots\x|\Fm_{q_l}^c|\x|\Fm_{r_1}^i|\x\dots\x|\Fm_{r_m}^i|$ by 
$|\Fm_{\vec q,\vec r}|$, where $\vec q=(q_1,\dots,q_l)$ and $\vec r=(r_1,\dots,r_m)$.

Given a $\lambda$-closed atomic $(n,\vec q,\vec r)$-meta-formula $\G=\vec p,\vec\gamma\mapsto(\vec x\mapsto\Rt F)$, 
where $\vec q=(q_q,\dots,q_l)$ and $\vec r=(r_1,\dots,r_m)$, its alternative 
interpretation will be a map $\KK\G\KK:|\Fm_{\vec q,\vec r}|\to\Hom(|\T|^n,\Qm)$, defined as follows.
Let us start from the usual Tarski-style interpretation 
$|\G|\:\Hom(\D^{\vec q,\vec r},\O)\to\Hom(\D^n,\Qm)$.
Given an $n$-tuple $\vec T$ of semantic $N_i$-terms $T_i\:\D^{\#N_i}\to\D$, we can combine them into a map 
$T\:\D^{\#N}\to\D^n$, where $N=N_1\cup\dots\cup N_n$.
Given an $l$-tuple $\vec F$ of semantic $(L_j,M^c_j,M^i_j,q_j)$-c-formulas 
\[F_j\:\D^{\#L_j}\x\Hom(\D^{\#M^c_j,\#M^i_j},\O)\to\Hom(\D^{q_j},\O_c)\] 
and an $m$-tuple $\vec\Phi$ of semantic $(\tilde L_k,\tilde M^c_k,\tilde M^i_k,r_k)$-i-formulas
\[\Phi_k\:\D^{\#\tilde L_k}\x\Hom(\D^{\#\tilde M^c_k,\#\tilde M^i_k},\O)\to\Hom(\D^{r_k},\O_i),\] 
we can combine them into a map 
\[F\x\Phi\:\D^{\#L}\x\Hom(\D^{\#M^c,\#M^i},\O)\to\Hom(\D^{\vec q,\vec r},\O),\] 
where $L=L_1\cup\dots\cup L_l\cup\tilde L_1\cup\dots\tilde L_m$, 
$M^c=M^c_1\cup\dots\cup M^c_l\cup\tilde M^c_1\cup\dots\cup\tilde M^c_m$ and 
$M^i=M^i_1\cup\dots\cup M^i_l\cup\tilde M^i_1\cup\dots\cup\tilde M^i_m$.
Finally, $|\G|$, $T$ and $F\x\Phi$ combine into a map $\D^{\#(L\cup N)}\x\Hom(\D^{\#M^c,\#M^i},\O)\to\Qm$, 
which we will denote, somewhat loosely, by $|\G|\circ(\vec T\x\vec F\x\vec\Phi)$.
We can apply to this map the Tarski-interpreted $\#(L\cup N)$-fold first-order and 
$(\#M^c,\#M^i)$-fold second-order meta-quantifier 
\[|\q_{\#(L\cup N)}^{\#M^c,\#M^i}|\:\Hom\big(\D^{\#(L\cup N)}\x\Hom(\D^{\#M^c,\#M^i},\O),\Qm\big)\to\Qm\]
and obtain an element $|\q_{\#(L\cup N)}^{\#M^c,\#M^i}|\big(|\G|\circ(\vec T\x\vec F\x\vec\Phi)\big)\in\Qm$.
Thus we have described a map $\KK\G\KK:|\Fm_{\vec q,\vec r}|\to\Hom(|\T|^n,\Qm)$.
Let us note that the Tarski-interpreted meta-quantifiers $|\q|$, $|\q^n_c|$ and $|\q^n_i|$ are used here
to interpret $\Rt$ rather than any meta-quantifiers.

Next, we keep the Tarski-style interpretations $\KK\&\KK:=|\&|$ and $\KK\impord\KK:=|\impord|$ for 
the meta-connectives, whereas the first-order meta-quantifier $\q$ is re-interpreted 
by a function $\KK\q\KK:\Hom(|\T|,\Qm)\to\Qm$; and the $n$-ary second-order meta-quantifiers $\q^n_c$, $q^n_i$ 
are re-interpreted by functions $\KK\q^n_c\KK\:\Hom(|\Fm_n^c|,\Qm)\to\Qm$ and
$\KK\q^n_i\KK\:\Hom(|\Fm_n^i|,\Qm)\to\Qm$.

The above will be called the {\it alternative meta-interpretation} associated to the given 
Tarski-style meta-interpretation.
In the case of a Tarski-style interpretation (i.e.\ two-valued meta-interpretation), we can 
define an {\it alternative interpretation} by providing explicit definitions of the functions
$\KK\q\KK$, $\KK\q^n_c\KK$ and $\KK\q^n_i\KK$.
Namely, we set $\KK\q\KK(f)=\Top$ if and only $f(T)=\Top$ for all $T\in|\T|$, $\KK\q^n_c\KK(f)=\Top$
if and only $f(F)=\Top$ for all $F\in|\Fm_n^c|$ and $\KK\q^n_i\KK(f)=\Top$ if and only 
$f(\Phi)=\Top$ for all $\Phi\in|\Fm_n^i|$.

Given an alternative meta-interpretation, we can interpret any $\lambda$-closed 
$(n,\vec q,\vec r)$-meta-formula $\G$ by a function $\KK\G\KK\:|\Fm_{\vec q,\vec r}|\to\Hom(|\T|^n,\Qm)$ 
by a straightforward induction.
In particular, any $\lambda$-closed meta-formula $\F$ is interpreted by an element $\KK\F\KK\in\Qm$.

\subsection{Intended informal semantics} \label{natural}

This is based on the alternative interpretation of the meta-logic and combines the ``meta-clarified 
BHK interpretation'' of intuitionistic logic \cite{M0}*{\S\ref{int:extended+}} with the corresponding 
version of the ``Verificationist interpretation'' of classical logic 
\cite{M0}*{Remark \ref{int:verificationist3}}.

$\O_c$ is taken to be a class of propositions, containing a prescribed class of contentful (e.g.\ mathematical) 
primitive propositions, and $\O_i$ is taken to be a class of problems, containing 
a prescribed class of contenful (e.g.\ mathematical) primitive problems.
Composite propositions and composite problems are obtained inductively from the primitive ones 
by using contentual classical and intuitionistic connectives and $\D$-indexed quantifiers, as well as
the contentual conversion operators $|\oc|\:\O_c\to\O_i$, $|\oc|(P)=$ ``Prove $P$'', and 
$|\wn|\:\O_i\to\O_c$, $|\wn|(\Gamma)=$ ``$\Gamma$ has a solution''.
The connectives, quantifiers and conversion operators of QHC are interpreted straightforwardly by 
their contentual analogues.

$\Qm$ is taken to be a class of problems containing the class $\O_i$ of ``meta-primitive'' problems.
Thus the function $\ocf_i\:\O_i\to\Qm$ is the inclusion; and the function $\ocf_c\:\O_c\to\Qm$ is 
the composition $\O_c\xr{|\oc|}\O_i\xr{\ocf_i}\Qm$.
``Meta-composite'' problems are obtained inductively from the meta-primitive ones by using contentual 
intuitionistic connectives and quantifiers over arbitrary $\D$-indexed, $\Hom(\D^n,\O_i)$-indexed 
($n=0,1,\dots$), $\Hom(\D^n,\O_c)$-indexed ($n=0,1,\dots$), $|\T|$-indexed, 
$|\Fm^i_n|$-indexed ($n=0,1,2,\dots$) and $|\Fm^c_n|$-indexed ($n=0,1,2,\dots$) families of problems.
Meta-connectives and meta-quantifiers are Tarski-interpreted by contentual intuitionistic connectives and 
$\D$-indexed, $\Hom(\D^n,\O_i)$-indexed and $\Hom(\D^n,\O_c)$-indexed quantifiers 
in the straightforward way; also, meta-quantifiers are alternatively interpreted by contentual intuitionistic 
$|\T|$-indexed, $|\Fm^i_n|$-indexed and $|\Fm^c_n|$-indexed quantifiers.

The function $\wnf\:\Qm\to\{\Top,\Bot\}$ sends a problem $\Gamma$ to $\Top$ if and only if there exists 
a solution of $\Gamma$.
What is meant by ``solving'' a problem $\Gamma$ is defined explicitly for primitive $\Gamma\in\O_i$;
also, ``truth'' is defined explicitly for primitive propositions $P\in\O_c$.
Then these notions are extended inductively to composite $P\in\O_c$ and $\Gamma\in\O_i$ by means of the usual 
truth tables and the usual BHK interpretation, along with two additional clauses: $|\wn|(\Gamma)$ is set to be 
``true'' if and only if $\Gamma$ has a solution; and to ``solve'' $|\oc|(P)$ means to prove that $P$ is true.
Finally, what is meant by ``solving'' a problem is further explained for all meta-composite $\Gamma\in\Qm$ by 
the usual BHK-interpretation.
This determines a contentual interpretation of QHC and its meta-logic, which we will call the {\it intended} 
interpretation.

\subsection{Modified informal semantics}

As noted above, Kolmogorov's interpretations of principles and rules are not compatible with the usual
Tarski-style theory.
But there is a simple workaround.
The following modified semantics overcomes this incompatibility by replacing the informal (i.e., contentual) 
constructive quantifier discussed in \S\ref{clarified} above with a combination of an informal classical 
quantifier (enabling compatibility with standard model theory) and a potentially formal constructive quantifier 
(preventing intuitionistic logic from degenerating into classical).
It is a straightforward combination of the ``Verificationist'' version (see \cite{M0}*{\S\ref{int:verificationist+}}) 
of the ``modified Platonist interpretation'' of classical logic \cite{M0}*{\S\ref{int:mPlatonist}} 
and the ``modified BHK interpretation'' of intuitionistic logic \cite{M0}*{\S\ref{int:modified+}}.

In addition to a domain of discourse $\D$, we fix a ``hidden parameter'' domain $\E$.
Now $\O_c$ is taken to be a class of unary predicates on $\E$, containing a prescribed class of contentful
(e.g.\ mathematical) primitive predicates, and $\O_i$ is taken to be a class of problems with parameter 
in $\E$, containing a prescribed class of contentful (e.g.\ mathematical) primitive parametric problems.
Composite predicates and composite parametric problems are obtained inductively from the primitive ones 
by using contentual classical and intuitionistic connectives and $\D$-indexed quantifiers, as well as
the contentual conversion operators $|\oc|\:\O_c\to\O_i$, $(|\oc|P)(e)=$ ``Prove $P(e)$'', and 
$|\wn|\:\O_i\to\O_c$, $(|\wn|\Gamma)(e)=$ ``$\Gamma(e)$ has a solution''.
The connectives, quantifiers and conversion operators of QHC are interpreted straightforwardly by 
their contentual analogues.

The value of $\ocf_i\:\O_i\to\Qm=\{\Top,\Bot\}$ on a parametric problem $\Gamma\in\O_i$ will be $\Top$ 
if and only if there exists a general method of solving the problem $\Gamma(e)$ for all $e\in\E$.
The function $\ocf_c\:\O_c\to\Qm$ is defined as the composition $\O_c\xr{|\oc|}\O_i\xr{\ocf_i}\Qm$.
(Thus the value of $\ocf_c$ on a predicate $P\in\O_c$ is $\Top$ if and only if {\it there exists a general 
method of proving that} $P(e)$ is true for all $e\in\E$, i.e., the passage from classical logic to classical 
meta-logic is through constructive logic.)

Here ``truth'' for $P(e)$ is defined explicitly for primitive $P\in\O_c$ and all $e\in\E$, and what is meant by 
``solving'' $\Gamma(e)$ is defined explicitly for primitive $\Gamma\in\O_i$ and all $e\in\E$ 
(separately for each $e\in\E$).
Then these notions are extended inductively to composite $P\in\O_c$ and $\Gamma\in\O_i$ (separately
for each $e\in E$) by means of the usual truth tables and the usual BHK interpretation, along with
two additional clauses: $(|\wn|\Gamma)(e)$ is set to be ``true'' if and only if $\Gamma(e)$ has a solution;
and to ``solve'' $(|\oc|P)(e)$ means to prove that $P(e)$ is true.
This determines a contentual interpretation of QHC and its meta-logic, which we will call 
the {\it modified} interpretation.

\section{Principles} \label{principles section}

\subsection{H- and K-Principles}\label{H-K}

\subsubsection{Hilbert's No Ignorabimus Principle}
Hilbert \cite{Hi1} (see also \cite{M0}*{\S\ref{int:controversy}}) prefaced his famous problem list
with an expression of his ``conviction'' that
\smallskip
\begin{center}
\parbox{14.7cm}{\small
 ``every definite mathematical problem must necessarily be susceptible of an
exact settlement, either in the form of an actual answer to the question asked,
or by the proof of the impossibility of its solution and therewith the necessary
failure of all attempts''.
}
\end{center}
\medskip
Thus in the terminology of \cite{M1}*{\ref{g1:semi-stable-decidable}}, he presumably asserted 
semi-decidability of all problems, given that he hardly intended to imply a general method \cite{Hi1}:
\smallskip
\begin{center}
\parbox{14.7cm}{\small
``This conviction of the decidability of every mathematical problem is a powerful incentive to
the worker. We hear within us the perpetual call: {\it There is the problem. Seek its solution.
You can find it by pure reason, for in mathematics there is no ignorabimus\footnotemark.}''
}
\footnotetext{Hilbert presumably refers to the Latin maxim {\it ignoramus et ignorabimus} (literally
``we do not know and shall not know'', in the sense ``[there are things] man can never know''),
and with it also likely to the then-famous position in the philosophy of science, articulated
by the physiologist Emil du Bois-Reymond (a brother of the mathematician Paul du Bois-Reymond)
in his 1872 address to the German Scientific Congress entitled ``\"Uber die Grenzen des Naturerkennens''
(``On the limits of our understanding of nature'') and in a number of later speeches and writings.}
\end{center}
\medskip
Hilbert reiterated his conviction 25 years later, with the ``perpetual call'' repeated verbatim, and
with the following addition \cite{Hi2}:

\smallskip
\begin{center}
\parbox{14.7cm}{\small
``Now, to be sure, my proof theory cannot specify a general method for deciding every mathematical
problem; that does not exist. But the demonstration\footnotemark\ that the assumption of the decidability of
every mathematical problem is consistent falls entirely within the scope of our theory.''
}
\footnotetext{This announced demonstration, which was to be ``the first and most important step'' in
Hilbert's attempted proof of the Continuum Hypothesis, seems to have never appeared.
Hilbert's paper contains a precise formulation of what he meant by ``decidability of every mathematical
problem'' only in ``the part that is of relevance here'', which concerned elimination of a choice
operator in formal proofs.}
\end{center}
\medskip

\noindent\formulas
In QHC, problems are represented by i-formulas, and so we will refer to the principle of decidability 
for i-formulas, $\prin\wn(\gamma\lor\neg\gamma)$ as ``Hilbert's No Ignorabimus Principle'', or 
the {\it H-principle} for brevity.
The principle of semi-decidability for i-formulas, $\prin\wn(\alpha\lor\neg\alpha)$,
is equivalent to the principle of semi-stability for i-formulas, 
$\prin\wn(\neg\neg\alpha\to\alpha)$, using that $\prin\neg\neg(\alpha\lor\neg\alpha)$
is derivable in intuitionistic logic (see \cite{M0}*{(\ref{int:not-not-LEM})}).
\metameta

A further argument in favour of the H-principle is that an i-formula of the form%
\footnote{In fact, there is no loss of generality here in considering only i-formulas of this form
(see \ref{Hilbert} below).}
$\oc F$ is semi-decidable as long as the c-formula $F$ satisfies either of Heyting's ``definitive''
judgement combinations (see \cite{M1}*{\S\ref{g1:letters1}}):
\begin{enumerate}[label=(\roman*)]
\item $\turnstile F$ (or, equivalently, $\turnstile\wn\oc F$);
\item $\turnstile\neg F$ (or, equivalently, $\turnstile\wn\oc\neg F$);
\item $\turnstile\neg\Box F$ and $\turnstile\neg\Box\neg F$ (or, equivalently,
$\turnstile\wn\neg\oc F$ and $\turnstile\wn\neg\oc\neg F$);
\item $\turnstile\neg\Box\neg\Box F$ and $\turnstile\neg\Box\neg\Box\neg F$  (or, equivalently,
$\turnstile\wn\neg\neg\oc F$ and $\turnstile\wn\neg\neg\oc\neg F$).
\end{enumerate}
Indeed, $\turnstile\wn\oc F\to\wn(\oc F\lor\neg\oc F)$, and also
$\turnstile\wn\oc\neg F\To\wn\neg\neg\oc\neg F$
and
$\turnstile\wn\neg\neg\oc\neg F\To\wn\neg\oc F$
and
$\turnstile\wn\neg\oc F\To\wn(\oc F\lor\neg\oc F)$.

An argument against the H-principle is that while it is quite reasonable to expect that any given
problem either has a solution, or can be proved (in some sense) to have no solutions --- there is no reason
to expect that the proof in question will always be expressible in a given formal system.
Yet this might well be required by the H-principle, depending on how one interprets $\oc$ (which is
implicit in the intuitionistic $\neg$).
(From this viewpoint, Heyting's ``definitive'' combinations are not really any more definitive than
the other ones, including the empty combination.)

\formulas
It should be noted also that the $\Diamond$-translation (see \cite{M1}*{\S\ref{g1:diamond}})
of $\prin\wn(\alpha\lor\neg\alpha)$ is $\prin\neg\wn\neg(\alpha\lor\neg\alpha)$, which is derivable in QHC
(since so is $\prin\neg\neg(\alpha\lor\neg\alpha)$).
\metameta

\subsubsection{Kolmogorov's Stability Principle}
Often a mathematical problem is formulated, either explicitly or implicitly, as {\it Prove or disprove $P$},
where $P$ is a proposition --- which amounts to $\oc P\lor\oc\neg P$ in our notation.
This is quite different from Hilbert's problem of ``exact settlement'', which in
the case of a problem of the form $\Gamma=\oc P$ specializes to $\oc P\lor\neg\oc P$.
In fact, any independent statement yields an example of the difference
(see \cite{M0}*{\S\ref{int:BHK-to}}).

However, in some cases of interest there is no difference:

\begin{example} Let $Q(n)$ denote the arithmetical predicate ``$2n$ is a sum of two primes''.
(Thus $\forall n\,Q(n)$ is Goldbach's conjecture, the example whose discussion by Heyting
and Kolmogorov was reviewed in \cite{M1}*{\S\ref{g1:letters2}}.)
Then for $P:=\neg\forall n\,Q(n)$, the negation of Goldbach's conjecture, there is 
no difference between $\oc\neg P$ and $\neg\oc P$.

Indeed, by checking all primes $\le 2n$ one can either prove or disprove $Q(n)$; thus $Q(n)$
is decidable: $\turnstile\forall n\big(\oc Q(n)\lor\oc\neg Q(n)\big)$.
Now $\oc\neg P$ is the problem $\oc\forall n\,Q(n)$ of proving Goldbach's conjecture;
if we restate it equivalently as $\forall n\,\oc Q(n)$, it requests us to find a general method of proving
$Q(n)$ for each $n$.
We actually do have such a general method (namely, checking all primes $\le 2n$), except that we do not
really know if it always works.
But it certainly cannot fail if we assume $\neg\oc P$, that is, $\neg\oc\neg\forall n\,Q(n)$.
In more detail, the latter problem is equivalent to $\neg\oc\exists n\,\neg Q(n)$, which by ($\octo$)
implies $\neg\exists n\,\oc\neg Q(n)$, which is in turn equivalent (see
\cite{M1}*{(\ref{int:quantifiers1})}) to $\forall n\,\neg\oc\neg Q(n)$.
But given that our checking procedure also solves the problem $\forall n\,\big(\oc Q(n)\lor\oc\neg Q(n)\big)$,
from this we get a solution of $\forall n\,\oc Q(n)$ --- the desired proof of Goldbach's conjecture.
\end{example}

What this argument really shows is that not only Goldbach's conjecture, but {\it every} arithmetical
$\Pi^0_1$ proposition%
\footnote{For instance, the usual formulation of Fermat's last theorem is another example of
an arithmetical $\Pi^0_1$ proposition. Riemann's hypothesis and the Poincare conjecture are known
to be equivalent to arithmetical $\Pi^0_1$ propositions by nontrivial but elementary results.}
(and still more generally, every universally quantified decidable proposition)
is stable --- that is, for its negation $P$ there is no difference between $\oc\neg P$ and $\neg\oc P$
(see \cite{M1}*{\S\ref{g1:stable and decidable}} concerning stable propositions).

\begin{remark}
Kent found an arithmetical formula $F$ such that in Peano Arithmetic,
$\turnstile F\to\Bew(F)$,
where $\Bew$ is G\"odel's provability predicate, but $F$ is not equivalent in Peano Aritmetic
to any $\Sigma^0_1$ formula (see \cite{Sm}*{4.3.2}).
If we let $\Box G:=G\land \Bew(\#G)$ (see \cite{M1}*{\S\ref{int:formal provability}}), 
then $\turnstile F\to\Box F$, or equivalently $\turnstile\Diamond\neg F\to\neg F$, which in turn implies
$\turnstile\Box\Diamond\neg F\to\neg F$.
(The latter is equivalent to the stability of $\neg F$, see \cite{M1}*{\S\ref{g1:stable and decidable}}.)
Thus in a certain arithmetical model of QS4, not all stable propositions are $\Pi^0_1$ propositions.
Beware that QS4 is not complete with respect to this model, which satisfies the independent principle
Grz of QS4 (see \cite{M0}*{\S\ref{int:formal provability}}).
\end{remark}

\formulas
As discussed in \cite{M1}*{\S\ref{g1:letters2}}, an opinion expressed by Kolmogorov in his letter to
Heyting may be interpreted as asserting that in a constructive framework, every sentence must be
either a problem or a stable proposition.
Since in QHC, propositions are represented by c-formulas,
we will use the title of ``Kolmogorov's Stablity Principle'', or the {\it K-principle}, for
the principle of stability for c-formulas, $\prin\neg\oc\neg p\to\oc p$,
--- or, equivalently (see \cite{M1}*{\ref{g1:semi-stable-decidable}(a)}),
the principle of semi-stability for c-formulas, $\prin\wn(\neg\oc\neg p\to\oc p)$.
\metameta

\subsubsection{Equivalent forms}
\formulas

We will use the following notation from \cite{M1}*{\S\ref{g1:stable and decidable}}:
$\Dec,\Stab:\1_i\too\1_i$ are defined by $\alpha\mapsto\alpha\lor\neg\alpha$ and 
$\alpha\mapsto\neg\neg\alpha\to\alpha$ respectively;
and $\Dec,\Stab:\1_c\too\1_i$ are defined by $p\mapsto\oc p\lor\oc\neg p$ and 
$p\mapsto\neg\oc\neg p\to\oc p$ respectively.

\begin{proposition} (a) The K-principle is equivalent to $\prin\oc\neg p\tofrom\neg\oc p$.

(b) The H-principle is equivalent to $\prin\wn\neg\alpha\tofrom\neg\wn\alpha$.
\end{proposition}

This follows directly from \cite{M1}*{\ref{g1:negation-commutes}}.

\begin{proposition} \label{Kolmogorov implies Hilbert}
(a) The H-principle is equivalent to $\prin\nabla\alpha\tofrom\neg\neg\alpha$.

(b) The K-principle implies $\prin\nabla\alpha\tofrom\neg\neg\alpha$.
\end{proposition}

\begin{proof} By the classical law of double negation,
$\turnstile\nabla\alpha\tofrom\oc\neg\neg\wn\alpha$.
Assuming the K-principle, we have
$\prin\oc\neg\neg\wn\alpha\tofrom\neg\neg\nabla\alpha$;
and assuming the H-principle, we have
$\prin\oc\neg\neg\wn\alpha\tofrom\nabla\neg\neg\alpha$.
In either case, by \cite{M1}*{\ref{g1:move-nabla}},
$\turnstile\nabla\neg\neg\alpha\tofrom\neg\neg\alpha$
and $\turnstile\neg\neg\alpha\tofrom\neg\neg\nabla\alpha$.
Thus in both cases we obtain $\prin\nabla\alpha\tofrom\neg\neg\alpha$.

It remains to prove that $\prin\nabla\alpha\tofrom\neg\neg\alpha$ implies
the H-principle: $\prin\nabla(\alpha\lor\neg\alpha)$.
But on assuming $\prin\nabla\beta\tofrom\neg\neg\beta$, the latter gets identified with the 
intuitionistically derivable principle $\prin\neg\neg(\alpha\lor\neg\alpha)$ 
(see \cite{M0}*{(\ref{int:not-not-LEM})}).
\end{proof}

\begin{proposition}\label{Kolmogorov} The following principles and rules are equivalent:

(i) the K-principle;

(ii) $\dfrac{\Diamond p}{p}$;

(iii) $\prin\Box p\tofrom p$;

(iv) $\prin\Box(p\lor q)\Tofrom(\Box p\lor\Box q)$;

(v) $\prin\wn\Dec(p)$.
\end{proposition}

As usual, $\Diamond p$ is an abbreviation for $\neg\Box\neg p$.

\begin{proof} We will establish the implications (i)$\imp$(ii)$\imp$(iii)$\imp$(iv)$\imp$(v)$\imp$(i).

From (i) we easily have $\prin\wn\neg\oc\neg p\to\wn\oc p$.
On the other hand, using that the K-principle implies the H-principle,
we can permute $\wn$ and $\neg$.
Thus we obtain $\prin\neg\wn\oc\neg p\to\wn\oc p$, or $\prin\Diamond p\to\Box p$.
Since $\turnstile\Box p\to p$, this in particular implies the rule (ii).

Next, it is well-known and easy to check the following in QS4 (and hence in QHC):
\begin{enumerate}
\item $\turnstile p\to\Diamond p$;
\item $\turnstile\Diamond\Diamond p\to\Diamond p$;
\item $\turnstile (\Diamond p\to\Diamond q)\to\Diamond(p\to q)$.
\end{enumerate}
From (3) we have $\turnstile (\Diamond\Diamond q\to\Diamond q)\to\Diamond(\Diamond q\to q)$.
Consequently, from (2) and {\it modus ponens} we get $\turnstile\Diamond(\Diamond q\to q)$.
Assuming the rule (ii), we infer from this $\prin\Diamond q\to q$.
This amounts to $\prin p\to\Box p$, by substituting $\neg p$ for $q$ (and $\neg q$ for $p$).

The principle (iv) trivially follows from (iii).

The law of excluded middle for c-formulas entails $\turnstile\Box (p\lor\neg p)$, which together with
the principle (iv) yields $\prin\Box p\lor\Box\neg p$.
Then from ($\wn_\lor$) we get $\prin\wn(\oc p\lor\oc\neg p)$; that is, (v).

Finally, semi-decidability implies semi-stability using that
$\turnstile\alpha\lor\beta\to\neg\beta\to\alpha$.
\end{proof}

\begin{remark} (a) It is also easy to get (ii) from (iii):
If $\prin p\to\Box p$, then $\prin \Diamond q\to q$, by substituting $\neg q$ for $p$.

(b) Here is a direct derivation of (v) from (i).
Assuming the K-principle, by the above we have $\prin\oc\neg p\tofrom\neg\oc p$ and the H-principle:
$\prin\wn(\alpha\lor\neg\alpha)$.
Thus we get $\prin\wn(\oc p\lor\neg\oc p)$, and consequently $\prin\wn(\oc p\lor\oc\neg p)$.
\end{remark}

\begin{proposition}\label{Hilbert} The following are equivalent:
\begin{enumerate}
\item the H-principle;
\item $\prin\wn\Dec(\nabla\alpha)$;
\item $\prin\Stab(\nabla\alpha)$;
\item $\prin\wn\Stab(\nabla\alpha)$.
\end{enumerate}
\end{proposition}

Here ``$\nabla\alpha$'' can clearly be replaced by ``$\oc p$''.
By \cite{M1}*{\ref{g1:stable-decidable}} and \cite{M1}*{\ref{g1:stable-decidable-bis}},
it can also be replaced by ``$\wn\alpha$''.
This can be in turn replaced by ``$\Box p$''.

\begin{proof}
(4) is trivially equivalent to $\prin\nabla(\neg\neg\nabla\alpha\to\nabla\alpha)$, which by
\cite{M1}*{\ref{g1:move-oc-wn2}(a)} is equivalent to
$\prin\neg\neg\nabla\alpha\to\nabla\alpha$, i.e., (3).
Now by \cite{M1}*{\ref{g1:nabla-negneg1}}, (3) is equivalent to
$\prin\neg\neg\alpha\to\nabla\alpha$.
By \ref{Kolmogorov implies Hilbert}(a), the latter is in turn equivalent to (1).

On the pother hand, by \cite{M1}*{\ref{g1:stable-decidable-bis}}, (2) is equivalent to 
$\prin\wn\Dec(\wn\alpha)$.
By \cite{M1}*{\ref{g1:semi-stable-decidable}(b)}, the latter is equivalent $\prin\wn\Dec(\alpha)$, 
i.e., (1).
\end{proof}

\begin{proposition}\label{diamond} The following are equivalent forms of the H-principle:
\begin{enumerate}
\item $\prin\Box p\tofrom\Diamond\Box p$;
\item $\prin\Box p\tofrom\Box\Diamond\Box p$.
\end{enumerate}
\end{proposition}

As noted in \cite{M1}*{\S\ref{g1:stable and decidable}}, the K-principle is equivalent to
$\prin\Box\Diamond p\to p$; by the above, it is also equivalent to $\prin p\tofrom\Box p$.

\begin{proof} The principle (1) is clearly equivalent to $\prin\wn\alpha\tofrom\Diamond\wn\alpha$, or
$\prin\wn\alpha\tofrom\neg\wn\oc\neg\wn\alpha$.
This implies $\prin\oc\wn\alpha\tofrom\oc\neg\wn\oc\neg\wn\alpha$; in fact, this implication is
reversible, by applying $\wn$ to both sides.
Thus (1) is equivalent to $\prin\nabla\alpha\tofrom\neg\neg\alpha$, a form of the H-principle.

By \ref{Hilbert} and the subsequent remark, the H-principle is equivalent to 
$\prin\Stab(\Box p)$; that is, to $\prin\neg\oc\neg\Box p\to\oc\Box p$.
By \cite{M1}*{\ref{g1:insolubility}}, this is equivalent to $\prin\oc\neg\wn\oc\neg\Box p\to\oc p$.
The latter implies $\prin\wn\oc\neg\wn\oc\neg\Box p\to\wn\oc p$; in fact, this implication is
reversible, by applying $\oc$ to both sides.
Thus the H-principle is equivalent to $\prin\Box\Diamond\Box p\to\Box p$.
The converse to the latter principle holds in QHC (in fact, already in QS4):
$\turnstile\Box\Box p\to\Box\Diamond\Box p$ due to $\turnstile q\to\Diamond q$.
\end{proof}

Let us note that one of the latest forms of the H-principle, $\prin\Diamond\Box p\to\Box p$, is equivalent to
$\prin\Diamond q\tofrom\Box\Diamond q$ by substituting $\neg q$ for $p$ (or $\neg p$ for $q$) and taking 
the contrapositive.
The latter is a law of the modal logic QS5, whose other laws can be taken to be 
those of QS4.
This yields

\begin{corollary} QHC extended by the H-principle is a conservative extension of QS5.
\end{corollary}

\begin{proof} Since the $\Box$-interpretation of QHC in QS4 is the identity on QS4, it descends
to an interpretation of QHC extended by the H-principle in QS5 that is the identity on QS5.
\end{proof}

\begin{proposition} \label{distrHK}
(a) The H-principle is equivalent to $\prin\wn(\alpha\to\beta)\Tofrom(\wn\alpha\to\wn\beta)$.

(b) The K-principle is equivalent to $\prin\oc(p\to q)\Tofrom(\oc p\to\oc q)$.
\end{proposition}

\begin{proof}[(a)] It suffices to show that $\prin\wn\neg\alpha\tofrom\neg\wn\alpha$ is not only
a special case of $\prin\wn(\alpha\to\beta)\tofrom(\wn\alpha\to\wn\beta)$, but also implies it.

Indeed, by the classical double negation law, $\wn(\alpha\to\beta)$ rewrites as $\neg\neg\wn(\alpha\to\beta)$;
and assuming that $\wn$ commutes with $\neg$, we can further rewrite it as $\neg\wn\neg(\alpha\to\beta)$.
By an intuitionistic validity (see \cite{M0}*{(\ref{int:implication1})}), the latter can be rewritten as
$\neg\wn(\neg\neg\alpha\land\neg\beta)$; and by ($\wn_\land$), the latter amounts to
$\neg(\wn\neg\neg\alpha\land\wn\neg\beta)$.
Once again using our assumption, we can rewrite the latter as $\neg(\neg\neg\wn\alpha\land\neg\wn\beta)$.
Validities of classical logic simplify this into $\wn\alpha\to\wn\beta$, as desired.
\end{proof}

\begin{proof}[(b)]
Assuming $\prin p\to\Box p$, we get $\prin\oc(p\to q)\tofrom(\oc p\to\oc q)$ from 
\cite{M1}*{\ref{g1:move-oc-wn}(a)}.
Conversely, $\prin\oc(p\to q)\tofrom(\oc p\to\oc q)$ specializes to $\prin\oc\neg p\tofrom\neg\oc p$.
\end{proof}

\metameta

 \subsection{Extraction principles}

\formulas

\begin{itemize}
\item {\it Proof Constructivity Principle} (PC-principle): $\prin\nabla\alpha\to\alpha$.

Informally, there is a method to extract from a proof of solubility of some problem an actual
solution of that problem.
\smallskip

\item {\it Reducibility Principle} (R-principle):
$\prin(\nabla\alpha\to\nabla\beta)\To \nabla(\alpha\to\beta)$.

Informally, there is a method to extract from a reduction of a proof of solubility of one problem to
a proof of solubility of another problem a proof of an actual reducibility of the first problem to
the second one.
\smallskip

\item {\it Disambiguation Principle} (D-principle):
$\prin\nabla(\alpha\lor\beta)\To (\nabla\alpha\lor\nabla\beta)$.

Informally, there is a method to extract from a proof of solubility of at least one of two problems
an explicit choice of one of these two problems whose solubility it actually proves.
\smallskip

\item {\it Proof Constructivity Rule} (PC-rule): $\dfrac{\nabla\alpha}\alpha$.
\smallskip

Informally, there is a method to find, for each formula $\Phi$, a method to extract from a proof of solubility 
(by a general method) of all problems instantiating $\Phi$ an actual solution (by a general method) of all 
those problems.
\end{itemize}
\metameta

\begin{remark}\label{refuting Hilbert} Let us note that if either $\turnstile\wn\Phi$ or
$\turnstile\neg\wn\Phi$, then $\turnstile\wn(\nabla\Phi\to\Phi)$.
(Indeed, by \cite{M1}*{\ref{g1:weak Hilbert} and \ref{g1:move-nabla}}, $\turnstile\neg\wn\Phi$ implies
$\turnstile\wn\neg\nabla\Phi$.)
Similarly, if $\Phi$ is an i-formula such that $\wn\Phi$ is either provable or disprovable in some theory
$T$ over QHC, then $T$ proves $\wn(\nabla\Phi\to\Phi)$.
\end{remark}

\formulas
\begin{proposition} \label{R+PC}
The meta-conjunction of the R-principle and the PC-rule is equivalent to the PC-principle.
\end{proposition}

\begin{proof} Clearly, $\prin\alpha\tofrom\nabla\alpha$ implies the R-principle and PC-rule.

Conversely, $\prin\nabla(\alpha\to\beta)\tofrom(\nabla\alpha\to\nabla\beta)$ implies, by
substituting $\nabla\beta$ for $\alpha$, the principle $\prin\nabla(\nabla\beta\to\beta)$, which
on applying the PC-rule yields $\prin\nabla\beta\to\beta$.
\end{proof}

\begin{proposition}\label{HK-nabla-box}
(a) The K-principle is equivalent to $\prin\Box(p\to q)\Tofrom(\Box p\to\Box q)$.

(b) The H-principle implies the R-principle, $\prin(\nabla\alpha\to\nabla\beta)\To \nabla(\alpha\to\beta)$.

\end{proposition}

\begin{proof}[Proof. (a)] $\prin p\tofrom\Box p$ trivially implies $\prin\Box(p\to q)\tofrom(\Box p\to\Box q)$.

Conversely, $\prin\Box(p\to q)\tofrom(\Box p\to\Box q)$ implies
$\prin(\Box p\to\Box q)\to(p\to q)$.
Substituting $\Box p$ for $q$ in the latter, we get $\prin p\to\Box p$.
\end{proof}

\begin{proof}[(b)] By applying $\oc$ to both sides of
$\prin\wn(\alpha\to\beta)\tofrom(\wn\alpha\to\wn\beta)$,
we get $\prin\nabla(\alpha\to\beta)\tofrom\oc(\wn\alpha\to\wn\beta)$, and the assertion follows from
\cite{M1}*{\ref{g1:move-oc-wn}(a)}.
\end{proof}

\begin{proposition} (a) The PC-rule, $\dfrac{\nabla\alpha}{\alpha}$, is equivalent to $\dfrac{\wn\alpha}{\alpha}$.\qquad

(b) The K-principle, $\dfrac{\Diamond p}{p}$, is equivalent to $\dfrac{\neg\oc p}{\neg p}$.
\end{proposition}

It should be noted that the second rule in (a) is the converse to ($\wn_\top$); and the second rule in (b)
is the converse to the rule form of ($\oc_\bot$) in \cite{M1}*{\ref{g1:v*}(b)}.

\begin{proof}
This is due to the reversibility of $(\oc_\top)$ and of the rule form of $(\wn_\bot)$ in
\cite{M1}*{\ref{g1:v*}(a)}.
\end{proof}

\begin{proposition} \label{principles}\

(a) The following are equivalent forms of the R-principle:
\begin{enumerate}[label=\oldstylenums{\arabic*}.,start=0]
\item $\prin\nabla(\alpha\to\beta)\Tofrom(\nabla\alpha\to\nabla\beta)$
\item $\prin\oc\wn(\alpha\to\beta)\Tofrom\oc(\wn\alpha\to\wn\beta)$;
\item $\prin\wn(\alpha\to\beta)\Tofrom\Box(\wn\alpha\to\wn\beta)$;
\item $\prin\wn(\alpha\to\beta)\Tofrom\wn(\nabla\alpha\to\nabla\beta)$;
\item $\prin\nabla(\alpha\to\beta)\Tofrom\nabla(\nabla\alpha\to\nabla\beta)$.
\end{enumerate}

(b) The following are equivalent forms of the D-principle:
\begin{enumerate}[label=\oldstylenums{\arabic*}.,start=0]
\item $\prin\nabla\alpha\lor\nabla\beta\Tofrom\nabla(\alpha\lor\beta)$
\item $\prin\oc\wn\alpha\lor\oc\wn\beta\Tofrom \oc(\wn\alpha\lor \wn\beta)$;
\item $\prin\oc p\lor\oc q\Tofrom \oc(\Box p\lor\Box q)$;
\item $\prin\oc p\lor \oc q\Tofrom\nabla(\oc p\lor \oc q)$;
\item $\prin\nabla\alpha\lor\nabla\beta\Tofrom\nabla(\nabla\alpha\lor\nabla\beta)$.
\end{enumerate}
\end{proposition}

It is easy to see that the ``$\to$'' implication in each line is derivable in QHC.
In fact, the ``$\to$'' in (a), \o2 is equivalent to ($\wn_\to$) by
using \cite{M1}*{\ref{g1:move-oc-wn2}(a)} (to see that it implies ($\wn_\to$), it suffices to use ($\wn\oc$));
and the ``$\to$'' in (b), \o2 is equivalent to ($\oc_\lor$) by using ($\wn\oc$).
	
\begin{proof}[Proof. (a)]
By \cite{M1}*{\ref{g1:move-oc-wn}(a)}, the right hand sides of \o0 and \o1 are equivalent
to each other; and the right hand sides of \o2 and \o3 are equivalent to each other.
By applying $\wn$ to both sides of \o1, we obtain \o2; conversely, by applying $\oc$ to both sides of \o2,
we recover \o1.
Similarly, by applying $\oc$ to both sides of \o3, we obtain \o4; and by applying $\wn$ to both sides
of \o4, we recover \o3.
(Alternatively, the right hand side of \o4 is equivalent to that of \o0 by
\cite{M1}*{\ref{g1:McKinsey+Gentzen}(e)}.)
\end{proof}

\begin{proof}[(b)]
By ($\wn_\lor$), the right hand sides of \o0 and \o1 are equivalent
to each other; and the right hand sides of \o2 and \o3 are equivalent to each other.
If we specialize \o1 or \o4 by substituting $\oc p$ for $\alpha$ and $\oc q$ for $\beta$, we obtain \o2 or
\o3, respectively; conversely, if we specialize \o2 or \o3 by substituting $\wn\alpha$ for $p$ and 
$\wn\beta$ for $q$, we recover \o1 or \o4, respectively.
(Alternatively, the right hand side of \o4 is equivalent to that of \o0 by
\cite{M1}*{\ref{g1:Goedel+Kuroda}(b)}.)
\end{proof}

\begin{corollary} \label{R-rule} The R-principle is equivalent to the rule 
$\dfrac{\wn\gamma\to(\wn\alpha\to\wn\beta)}{\wn\gamma\to\wn(\alpha\to\beta)}$.
\end{corollary}

A version of this rule appears in \cite{M0}*{\S\ref{int:weak BHK}} in connection with
the implication clause of the BHK interpretation.
Let us note that the said rule is not admissible in QHC since it is equivalent to a non-derivable (see
\S\ref{models} below) principle.

\begin{proof} The rule in question may be equivalently rewritten by applying $\oc$'s to both
the premise and the conclusion.
Using \cite{M1}*{\ref{g1:move-oc-wn}(a)} (the strengthened version which is actually proved),
this can further be equivalently rewritten as
\[\dfrac{\nabla\gamma\to(\nabla\alpha\to\nabla\beta)}{\nabla\gamma\to\nabla(\alpha\to\beta)}\text.\]
The latter rule clearly follows from the R-principle.
For the converse, let us substitute $\nabla\alpha\to\nabla\beta$ for $\gamma$; then the premise
is derivable by \cite{M1}*{\ref{g1:McKinsey+Gentzen}(e)} and the conclusion is equivalent to
the R-principle by \ref{principles}(a).
\end{proof}

\begin{proposition}\label{principles2}\

(a) The following are equivalent forms of the H-principle:
\begin{enumerate}[label=\oldstylenums{\arabic*}.]
\item $\prin(\wn\oc p\to \wn\oc q)\Tofrom\wn(\oc p\to\oc q)$;
\item $\prin(\wn\alpha\to\wn\beta)\Tofrom \wn(\nabla\alpha\to\nabla\beta)$;
\item $\prin(\wn\alpha\to\wn\beta)\Tofrom\Box(\wn\alpha\to\wn\beta)$;
\item $\prin(\Box p\to\Box q)\Tofrom\Box(\Box p\to\Box q)$.
\end{enumerate}

(b) The following are equivalent forms of the K-principle:
\begin{enumerate}[label=\oldstylenums{\arabic*}.,start=0]
\item $\prin\Box(p\lor q)\Tofrom\Box p\lor\Box q$;
\item $\prin\wn\oc(p\lor q)\Tofrom \wn(\oc p\lor \oc q)$;
\item $\prin\oc(p\lor q)\Tofrom\nabla(\oc p\lor \oc q)$;
\item $\prin\oc(p\lor q)\Tofrom\oc(\Box p\lor\Box q)$;
\item $\prin\Box(p\lor q)\Tofrom\Box(\Box p\lor\Box q)$.
\end{enumerate}

(c) The following are equivalent forms of the K-principle:
\begin{enumerate}[label=\oldstylenums{\arabic*}.,start=0]
\item $\prin\Box(p\to q)\Tofrom(\Box p\to\Box q)$;
\item $\prin\wn\oc(p\to q)\Tofrom\wn(\oc p\to\oc q)$;
\item $\prin\oc(p\to q)\Tofrom\nabla(\oc p\to \oc q)$;
\item $\prin\oc(p\to q)\Tofrom\oc(\Box p\to\Box q)$;
\item $\prin\Box(p\to q)\Tofrom\Box(\Box p\to\Box q)$.
\end{enumerate}
\end{proposition}

It is easy to see that the ``$\from$'' implication in each line of (a,b) and the ``$\to$'' implication in
each line of (c) is derivable in QHC.
In fact, the ``$\from$'' in (b), \o2 is equivalent to ($\oc_\lor$) by using
\cite{M1}*{\ref{g1:move-oc-wn2}(a)}.

The following ``missing'' analogue of \o2's is by \cite{M1}*{\ref{g1:move-oc-wn}(a)} derivable in QHC:
$$\turnstile(\oc p\to\oc q)\Tofrom \oc(\Box p\to\Box q).$$

Let us note that a special case of (a), \o4, $\prin\neg\Box p\tofrom\Box\neg\Box p$, is essentially the same
as one of the previous forms of the H-principle, $\prin\Box p\tofrom\Diamond\Box p$.

\begin{proof}[Proof. (a)] The equivalences \o1$\iff$\o2 and \o3$\iff$\o4 are similar to \ref{principles}(b).
Also, by \cite{M1}*{\ref{g1:move-oc-wn}(a)}, the right hand side of \o2 is equivalent to that of \o3.

Now \o3 specializes to $\prin\neg\wn\alpha\tofrom\Box\neg\wn\alpha$, whereas
$\turnstile \Box\neg\wn\alpha\tofrom\wn\neg\alpha$ due to \cite{M1}*{\ref{g1:insolubility}}.
Thus \o3 implies the H-principle.

Conversely, by \ref{HK-nabla-box}(a), the H-principle implies the R-principle, which by \ref{principles}(a)
is equivalent to $\prin\wn(\alpha\to\beta)\tofrom\wn(\nabla\alpha\to\nabla\beta)$.
On the other hand, by \ref{distrHK}(a), the H-principle also implies 
$\prin\wn(\alpha\to\beta)\tofrom(\wn\alpha\to\wn\beta)$.
By combining these, we get \o3.
\end{proof}

\begin{proof}[(b)]
This is proved similarly to \ref{principles}(a), using ($\wn_\lor$) to obtain the two equivalences
of the right hand sides.
\end{proof}

\begin{proof}[(c)]
The equivalence of \o1 and \o2 is similar to \ref{principles}(a).
The right hand side of \o3 is equivalent to that of \o2 by \cite{M1}*{\ref{g1:worst}(c)}, and to that of
$\prin\oc(p\to q)\Tofrom (\oc p\to\oc q)$ (yet another form of the K-principle)
by \cite{M1}*{\ref{g1:move-oc-wn}(a)}.
\end{proof}

\begin{proposition}\label{expressive} \

(a) The following are equivalent forms of the PC-principle, $\prin\nabla\alpha\tofrom\alpha$:
\begin{enumerate}[label=\oldstylenums{\arabic*}.]
\item $\prin(\alpha\to \beta)\Tofrom \oc(\wn\alpha\to\wn\beta)$;
\item $\prin\alpha\lor\beta\Tofrom\oc(\wn\alpha\lor \wn\beta)$,
\end{enumerate}

(b) The following are equivalent forms of the K-principle, $\prin\Box p\tofrom p$:
\begin{enumerate}[label=\oldstylenums{\arabic*}.]
\item $\prin(p\to q)\Tofrom \wn(\oc p\to\oc q)$;
\item $\prin p\lor q\Tofrom\wn(\oc p\lor \oc q)$.
\end{enumerate}
\end{proposition}

Here the ``$\to$'' implication in (a), \o2 and the ``$\from$'' implication in (b), \o2 are equivalent to
($\wn_\lor$) and ($\oc_\lor$) respectively by using \cite{M1}*{\ref{g1:Galois}}.
The ``$\to$'' implication in (a), \o1 is derivable in QHC by the same token.

\begin{proof}[Proof. (a)] By \cite{M1}*{\ref{g1:move-oc-wn}(a)}, \o1 is
equivalent to $\prin(\alpha\to\beta)\tofrom(\nabla\alpha\to\nabla\beta)$.
The latter principle clearly follows from $\prin\nabla\beta\tofrom\beta$; and also implies it by 
substituting $\triv$ for $\alpha$.

By ($\wn_\lor$), \o2 is equivalent to $\prin\alpha\lor\beta\tofrom\nabla(\alpha\lor\beta)$.
The latter principle clearly follows from $\prin\nabla\alpha\tofrom\alpha$; and also implies it by 
substituting either $\alpha$ or $\ab$ for $\beta$.
\end{proof}

\begin{proof}[(b)] Assuming $\prin p\tofrom\Box p$, \o1 follows from \ref{principles2}(c), \o1.
Conversely, by applying $\oc$ to both sides of \o1, we obtain \ref{principles2}(c), \o2.

By ($\wn_\lor$), \o2 is equivalent to $\prin p\lor q\tofrom(\Box p\lor\Box q)$.
This clearly follows from $\prin p\tofrom\Box p$; and also implies it by substituting either $p$
or $\clbot$ for $q$.
\end{proof}

\begin{proposition}
(a) QHC extended by any collection of consequences of the PC-principle, $\prin\nabla\alpha\tofrom\alpha$,
is conservative over QH and over QS4.

(b) QHC extended by any collection of consequences of the K-principle, $\prin\Box p\tofrom p$,
is conservative over QH and over QC.
\end{proposition}

\begin{proof} The PC-principle (and any of its consequences) is validated in the $\Box$-interpretation
of QHC in QS4, which is exact on both QH and QS4.
(Alternatively, we may use the $\nabla$-interpretation of QHC in itself.)

The K-principle (and any of its consequences) is validated in the $\neg\neg$-interpretation
of QHC in QH, which is exact on both QH and QC.
(Alternatively, we may use the $\Diamond$-interpretation of QHC in itself.)
\end{proof}

\metameta

\subsection{Principles involving quantifiers}

\formulas

\begin{proposition}\label{principles-q}\

(a) The following are equivalent forms of the PC-principle, $\prin\beta\tofrom\nabla\beta$:
\begin{enumerate}[label=\oldstylenums{\arabic*}.]
\item $\prin\exists\tr x\,\alpha(\tr x)\Tofrom \oc\exists\tr x\, \wn\alpha(\tr x)$;
\item $\prin\forall\tr x\, \alpha(\tr x)\Tofrom \oc\forall\tr x\, \wn\alpha(\tr x)$.
\end{enumerate}

(b) The following are equivalent forms of the K-principle, $\prin q\tofrom\Box q$:
\begin{enumerate}[label=\oldstylenums{\arabic*}.]
\item $\prin\exists\tr x\, p(\tr x)\Tofrom \wn\exists\tr x\, \oc p(\tr x)$;
\item $\prin\forall\tr x\, p(\tr x)\Tofrom \wn\forall\tr x\, \oc p(\tr x)$.
\end{enumerate}
\end{proposition}

Here the ``$\to$'' implications in (a) are equivalent to ($\wn_\exists$) and ($\wn_\forall$),
and the ``$\from$'' implications in (b) are equivalent to ($\oc_\exists$) and ($\oc_\forall$) respectively,
by using \cite{M1}*{\ref{g1:Galois}}.

\begin{proof} By using ($\wn_\exists$) and ($\oc_\forall$), we can rewrite each of the four principles
in the form where one side differs from another by a single $\Box$ or $\nabla$ (either before or after
the quantifier).
Hence the principles in (a) immediately follow from $\prin\beta\tofrom\nabla\beta$, and the principles in (b)
from $\prin q\tofrom\Box q$.
The converse implications are obtained by substituting $\tr x\mapsto\beta$ for $\alpha$, where $\beta$ is 
a nullary problem variable, and $\tr x\mapsto q$ for $p$, where $q$ is a propositional (i.e., nullary 
predicate) variable.
\end{proof}

\begin{proposition}\label{nabla-K}
 The following are equivalent to the K-principle:
\begin{enumerate}[label=\oldstylenums{\arabic*}.]
\item $\prin\exists\tr x\, p(\tr x)\Tofrom \neg\wn\neg\exists\tr x\, \oc p(\tr x)$;
\item $\prin\forall\tr x\, q(\tr x)\Tofrom \wn\forall\tr x\, \neg\oc\neg q(\tr x)$.
\end{enumerate}
\end{proposition}

This essentially asserts the equivalence between formulas of classical logic and their
$\Diamond$-translations (see \cite{M1}*{\S\ref{g1:diamond}}).

\begin{proof}
\o1 is equivalent to \o2 by substituting $\neg p$ for $q$ (or $\neg q$ for $p$), taking the contrapositive
and using \cite{M0}*{(\ref{int:quantifiers1})}.

It is easy to see that the K-principle implies \o2 (using that $\oc$ commutes with $\forall$ and,
on assuming the K-principle, with $\neg$, and that, on assuming the K-principle, $\Box$ can be erased).
Conversely, by substituting $\tr x\mapsto r$ for $q$, where $r$ is a propositional variable, from \o2 
we get $\prin r\tofrom\wn\neg\oc\neg r$.
Using that $\turnstile\neg\alpha\tofrom\oc\neg\wn\alpha$, we can rewrite this as 
$\prin r\tofrom\Box\Diamond r$.
Yet $\prin\Box\Diamond r\to r$ is one of the forms of the K-principle.
\end{proof}

The principles of the following proposition follow trivially from either
the K-principle, in the form $\prin p\tofrom\Box p$ (in the case of (b) and (c)), or from the PC-principle,
$\prin\alpha\tofrom\nabla\alpha$ (in the case of (a) and (d)).

\begin{proposition} \label{principles-q2}\

(a) The following are equivalent forms of the {\rm $\forall_\nabla$-principle},
$\prin\nabla\forall\tr x\, \alpha(\tr x)\Tofrom \forall\tr x\, \nabla \alpha(\tr x)$:
\begin{enumerate}[label=\oldstylenums{\arabic*}.]
\item $\prin\oc\wn\forall\tr x\,\alpha(\tr x)\Tofrom\oc\forall\tr x\,\wn\alpha(\tr x)$;
\item $\prin\wn\forall\tr x\,\alpha(\tr x)\Tofrom\Box\forall\tr x\, \wn\alpha(\tr x)$;
\item $\prin\wn\forall\tr x\, \alpha(\tr x)\Tofrom \wn\forall\tr x\, \nabla \alpha(\tr x)$;
\item $\prin\nabla\forall\tr x\, \alpha(\tr x)\Tofrom \nabla\forall\tr x\, \nabla \alpha(\tr x)$.
\end{enumerate}

(b) The following are equivalent forms of the {\rm $\exists_\Box$-principle},
$\prin\Box\exists\tr x\, p(\tr x)\Tofrom \exists\tr x\, \Box p(\tr x)$:
\begin{enumerate}[label=\oldstylenums{\arabic*}.]
\item $\prin\wn\oc\exists\tr x\, p(\tr x)\Tofrom \wn\exists\tr x\, \oc p(\tr x)$;
\item $\prin\oc\exists\tr x\, p(\tr x)\Tofrom\nabla\exists\tr x\, \oc p(\tr x)$;
\item $\prin\oc\exists\tr x\, p(\tr x)\Tofrom \oc\exists\tr x\, \Box p(\tr x)$;
\item $\prin\Box\exists\tr x\, p(\tr x)\Tofrom \Box\exists\tr x\, \Box p(\tr x)$.
\end{enumerate}

(c) The following are equivalent forms of the {\rm $\forall_\Box$-principle},
$\prin\forall\tr x\,\Box p(\tr x)\Tofrom \Box\forall\tr x\, p(\tr x)$:
\begin{enumerate}[label=\oldstylenums{\arabic*}.]
\item $\prin\forall\tr x\,\wn\oc p(\tr x)\Tofrom\wn\forall\tr x\,\oc p(\tr x)$;
\item $\prin\forall\tr x\,\wn\alpha(\tr x)\Tofrom \wn\forall\tr x\,\nabla\alpha(\tr x)$;
\item $\prin\forall\tr x\, \wn\alpha(\tr x)\Tofrom \Box\forall\tr x\, \wn\alpha(\tr x)$;
\item $\prin\forall\tr x\, \Box p(\tr x)\Tofrom \Box\forall\tr x\, \Box p(\tr x)$.
\end{enumerate}

(d) The following are equivalent forms of the {\rm $\exists_\nabla$-principle},
$\prin\exists\tr x\,\nabla\alpha(\tr x)\Tofrom\nabla\exists\tr x\,\alpha(\tr x)$:
\begin{enumerate}[label=\oldstylenums{\arabic*}.]
\item $\prin\exists\tr x\, \oc\wn\alpha(\tr x)\Tofrom \oc\exists\tr x\, \wn\alpha(\tr x)$;
\item $\prin\exists\tr x\,\oc p(\tr x)\Tofrom \oc\exists\tr x\,\Box p(\tr x)$;
\item $\prin\exists\tr x\, \oc p(\tr x)\Tofrom \nabla\exists\tr x\, \oc p(\tr x)$;
\item $\prin\exists\tr x\, \nabla\alpha(\tr x)\Tofrom \nabla\exists\tr x\, \nabla\alpha(\tr x)$.
\end{enumerate}
\end{proposition}

It is easy to see that the ``$\to$'' implication in each line of (a,d) and the ``$\from$''
implication in each line of (b,c) is derivable in QHC.
In fact, it is equivalent to the corresponding law of QHC, ($\wn_\forall$) or ($\wn_\exists$)
(by using \cite{M1}*{\ref{g1:move-oc-wn2}} in the cases of (a), \o2 and (b), \o2 and
using ($\oc\wn$) and ($\wn\oc$) in the cases of (c), \o2 and (d), \o2).

\begin{proof}[Proof. (a,b)]
These are proved similarly to \ref{principles}(a), using ($\oc_\forall$) and ($\wn_\exists$) respectively.
\end{proof}

\begin{proof}[(c,d)]
These are proved similarly to \ref{principles}(b), using ($\oc_\forall$) and ($\wn_\exists$) respectively.
\end{proof}

\begin{corollary} The $\forall_\nabla$-principle is equivalent to the rule 
$\dfrac{\wn\gamma\to\forall\tr x\,\wn\alpha(\tr x)}{\wn\gamma\to\wn\forall\tr x\,\alpha(\tr x)}$.
\end{corollary}

A version of this rule appears in \cite{M0}*{\S\ref{int:weak BHK}} in connection with
the universal quantifier clause of the BHK interpretation.
Let us note that the said rule is not admissible in QHC since it is equivalent to a non-derivable (see
\S\ref{models} below) principle.

The proof is similar to that of \ref{R-rule}.

Each of the $\forall_\Box$-, $\exists_\Box$-, $\forall_\nabla$- and $\exists_\nabla$-principles is
satisfied in every model of QHC with singleton domain; therefore none of these principles implies
any of the independent quantifier-free principles.

\begin{proposition} \label{h-quant}
The H-principle implies the $\forall_\Box$-principle.
\end{proposition}

\begin{proof} By the classical law of excluded middle,
$\turnstile\exists\tr x\,\Box r(\tr x)\lor\neg\exists\tr x\,\Box r(\tr x)$.
On the other hand, by \ref{principles2}(a), \o4, the H-principle implies $\prin\neg\Box q\tofrom\Box\neg\Box q$.
Since $\turnstile\exists\tr x\,\Box r(\tr x)\tofrom\Box\exists\tr x\,\Box r(\tr x)$ 
\cite{M1}*{\ref{g1:Goedel+Kuroda}(c)},
we get $\prin\exists\tr x\,\Box r(\tr x)\lor\Box\neg\exists\tr x\,\Box r(\tr x)$.
The latter principle is classically equivalent to 
$\prin\neg\forall\tr x\,\neg\Box r(\tr x)\lor\Box\forall\tr x\,\neg\Box r(\tr x)$,
and hence to $\prin\forall\tr x\,\neg\Box r(\tr x)\to\Box\forall\tr x\,\neg\Box r(\tr x)$.
The latter specializes, by substituting $\neg q(\tr x)$ for $r(\tr x)$, to 
$\prin\forall\tr x\,\Diamond q(\tr x)\to\Box\forall\tr x\,\Diamond q(\tr x)$,
which in turn specializes, by substituting $\Box p(\tr x)$ for $q(\tr x)$, to
$\prin\forall\tr x\,\Diamond\Box p(\tr x)\to\Box\forall\tr x\,\Diamond\Box p(\tr x)$.
By applying again the H-principle in the form of \ref{diamond}, $\prin\Diamond\Box q\tofrom\Box q$,
we get $\prin\forall\tr x\,\Box p(\tr x)\to\Box\forall\tr x\,\Box p(\tr x)$, 
i.e. the $\forall_\Box$-principle in the form of \ref{principles-q2}(c),\o4.
\end{proof}

\begin{proposition} \label{quant-quant}
The $\exists_\Box$-principle implies the $\forall_\Box$-principle.
\end{proposition}

\begin{proof} The $\exists_\Box$-principle implies, in particular, the rule
$$\frac{\exists\tr x\,q(\tr x)}{\exists\tr x\,\Box q(\tr x)}.$$
We will show that this rule implies the $\forall_\Box$-principle,
$\prin\forall\tr x\, \Box p(\tr x)\to\Box\forall\tr x\, \Box p(\tr x)$.

Let $\mm F(\tr x)=\neg\Box p(\tr x)\lor\forall\tr y\,\Box p(\tr y)$.
Then $\prin\exists\tr x\,\mm F(\tr x)$ is equivalent to 
$\prin\exists\tr x\,\neg\Box p(\tr x)\lor\forall\tr y\,\Box p(\tr y)$,
which is in turn equivalent to $\prin\neg\forall\tr x\,\Box p(\tr x)\lor\forall\tr x\,\Box p(\tr x)$ --- a special 
case of the classical law of excluded middle.
Thus $\prin\exists\tr x\,\mm F(\tr x)$ is derivable in QS4.
Our hypothesis then implies $\prin\exists\tr x\,\Box\mm F(\tr x)$, or in more detail,
$\prin\exists\tr x\,\Box\big(\neg\Box p(\tr x)\lor\forall\tr y\,\Box p(\tr y)\big)$.
This can be rewritten as $\prin\exists\tr x\,\Box\big(\Box p(\tr x)\to\forall\tr y\,\Box p(\tr y)\big)$,
which in turn implies $\prin\exists\tr x\,\big(\Box p(\tr x)\to\Box\forall\tr y\,\Box p(\tr y)\big)$.
By an intuitionistic validity (see \cite{M0}*{(\ref{int:gq2})}), the latter implies
$\prin\forall\tr x\,\Box p(\tr x)\to\Box\forall\tr y\,\Box p(\tr y)$.
This is the same as $\prin\forall\tr x\,\Box p(\tr x)\to\Box\forall\tr x\,\Box p(\tr x)$, 
i.e. the $\forall_\Box$-principle in the form of \ref{principles-q2}(c),\o4.
\end{proof}

\metameta

\begin{proposition} \label{quant-no quant}
The $\forall_\Box$-principle does not imply the $\exists_\Box$-principle.
\end{proposition}

\begin{proof} Under the $\Box$-interpretation, the two principles, when expressed using only
$\Box=\wn\oc$ (with no separate $\wn$'s and $\oc$'s), turn into principles of QS4 expressed by
the same formulas.
Every topological model of QS4 in an Alexandroff space satisfies the $\forall_\Box$-principle,
in the form $\fm{\prin\Box\forall\tr x\,\Box p(\tr x)\tofrom\forall\tr x\,\Box p(\tr x)}$.
The $\exists_\Box$-principle, in the form 
$\fm{\prin\Box\exists\tr x\,p(\tr x)\tofrom\Box\exists\tr x\,\Box p(\tr x)}$, fails,
for instance, in a topological model in the Alexandroff space corresponding to the poset $\N$ (with
the usual order), with $\D=\N$ and with $|\fm p|(n)=\{n\}$.
\end{proof}

\begin{proposition} \label{no-implic}
The $\forall_\nabla$-principle and the $\exists_\nabla$-principle do not imply one another.

Moreover, the D-principle implies neither the $\forall_\nabla$-principle nor the $\exists_\nabla$-principle.
\end{proposition}

\begin{proof} The $\neg\neg$-interpretation transforms the first two principles in question into
the $\neg\neg$-Shift Principle and the Strong Markov Principle (see \cite{M1}).
But these do not imply one another, as shown by Tarski models.
Namely, the $\neg\neg$-Shift Principle is satisfied in the one-point compactification $\N_+$ of
the countable discrete space, and fails in the Alexandroff space corresponding to the poset $\N$; whereas
the Strong Markov Principle is satisfied in all Alexandroff spaces whose corresponding preorder is
a directed set, and fails in $\N_+$, since the union of the regular open sets $\{n\}$, $n\in\N$, is
not regular open (see \cite{M0}*{\S\ref{int:Markov} and \S\ref{int:shift}}).

The second assertion follows since the $\neg\neg$-translation transforms the D-principle
into an equivalent form of de Morgan's principle, which does not imply the $\neg\neg$-Shift Principle
since it is also satisfied in Alexandroff spaces whose corresponding preorder is a directed set, and
does not imply the Strong Markov Principle since it is satisfied in the Stone--\v Cech compactification
of $\N$ (whose existence depends on the uncountable axiom of choice), or alternatively by analyzing
an obvious model in the one-point compactification of $\N$ (see \cite{M0}*{\S\ref{int:Markov} and
Proposition \ref{int:Stone-Cech}}).
\end{proof}

\begin{remark} An alternative way to see that the $\forall_\Box$- and $\forall_\nabla$-principles
do not imply the $\exists_\Box$- and $\exists_\nabla$-principles is by considering two-element domains,
over which the former principles hold, since $\land$ commutes with $\Box$ and $\nabla$, and the latter
principles generally do not hold, since $\lor$ generally does not commute with either $\Box$ or $\nabla$
(by Theorem \ref{intermediate models}).
\end{remark}

\subsection{Intuitionistically unacceptable principles}

\formulas
All independent principles of intuitionistic logic discussed in this section can be found
in \cite{M0}*{\S\ref{int:topology}}, which includes Tarski topological models that prove
their independence.

\begin{proposition} Assume the H-principle. Then:

(a) the PC-rule, $\dfrac{\nabla\alpha}{\alpha}$, is equivalent to the problem decidability principle
(i.e.\ $\prin\alpha\lor\neg\alpha$);

(b) the $\forall_\nabla$-principle is equivalent to the $\neg\neg$-Shift Principle;

(c) the $\exists_\nabla$-principle is equivalent to the Strong Markov Principle.
\end{proposition}

Since the K-principle is conservative over intuitionistic logic and implies the H-principle,
we obtain that the K-principle does not imply either the $\forall_\nabla$- or the $\exists_\nabla$-principle.

\begin{proof} Parts (b) and (c) follow since the H-principle is equivalent
to $\prin\nabla\alpha\tofrom\neg\neg\alpha$.
To get (a), we additionally use that the H-principle implies the principle
$\prin\nabla(\alpha\to\beta)\tofrom(\nabla\alpha\to\nabla\beta)$, whose meta-conjunction with
$\nabla\alpha/\alpha$ is equivalent to $\prin\nabla\alpha\tofrom\alpha$.
Thus on assuming the H-principle, $\nabla\alpha/\alpha$ is equivalent to $\prin\alpha\tofrom\neg\neg\alpha$,
which is in turn equivalent to $\prin\alpha\lor\neg\alpha$ (see \cite{M0}*{(\ref{int:decidable-stable}) and
(\ref{int:not-not-LEM})}).
\end{proof}

\begin{proposition} Assume the PC-principle. Then:

(a) the $\forall_\Box$-principle implies the Constant Domain Principle;

(b) the $\exists_\Box$-principle implies the Principle of Independence of Premise.
\end{proposition}

Since the PC-principle is conservative over intuitionistic logic, we obtain that the PC-principle
does not imply either the $\forall_\Box$- or the $\exists_\Box$-principle.

\begin{proof}[Proof. (a)]
By \ref{principles-q2}(c), the $\forall_\Box$-principle is equivalent to
$\prin\forall\tr x\,\wn\alpha(\tr x)\tofrom \wn\forall\tr x\,\nabla\alpha(\tr x)$, which in the presence of
$\prin\nabla\alpha\tofrom\alpha$ simplifies to
$\prin\forall\tr x\,\wn\alpha(\tr x)\tofrom \wn\forall\tr x\,\alpha(\tr x)$.

Since the Constant Domain Principle holds classically, we have
$$\turnstile \oc\forall\tr x\,(\wn\alpha\lor\wn\beta(\tr x))\Tofrom\oc(\wn\alpha\lor\forall\tr x\,\wn\beta(\tr x)).$$
Here the left hand side is equivalent by ($\oc_\forall$) and ($\wn_\lor$) to
$\forall\tr x\,\oc\wn(\alpha\lor\beta(\tr x))$, which in the presence of $\prin\nabla\alpha\tofrom\alpha$
simplifies to $\forall\tr x\,(\alpha\lor\beta(\tr x))$.
On the other hand, $\alpha\lor\forall\tr x\,\beta(\tr x)$ can be rewritten in the presence of
$\prin\nabla\alpha\tofrom\alpha$ as $\oc\wn(\alpha\lor\forall\tr x\,\beta(\tr x))$,
which by ($\wn_\lor$) is equivalent to $\oc(\wn\alpha\lor\wn\forall\tr x\,\beta(\tr x))$.
Collecting the previous steps, we get that the Constant Domain Principle,
$\prin\forall\tr x\,(\alpha\lor\beta(\tr x))\To\alpha\lor\forall\tr x\,\beta(\tr x)$, is equivalent
(in the presence of the PC-principle) to
$\prin\oc(\wn\alpha\lor\forall\tr x\,\wn\beta(\tr x))\To\oc(\wn\alpha\lor\wn\forall\tr x\,\beta(\tr x))$; and
the assertion follows.
\end{proof}

\begin{proof}[(b)]
By \ref{principles-q2}(b), the $\exists_\Box$-principle is equivalent to
$\prin\oc\exists\tr x\, p(\tr x)\tofrom\nabla\exists\tr x\, \oc p(\tr x)$, which in the presence of
$\prin\nabla\alpha\tofrom\alpha$ simplifies to $\prin\oc\exists\tr x\, p(\tr x)\tofrom\exists\tr x\, \oc p(\tr x)$.

Since the Principle of Independence of Premise holds classically, we have
$$\turnstile\oc(\wn\alpha\to\exists\tr x\,\wn\beta(\tr x))\Tofrom\oc\exists\tr x\,(\wn\alpha\to\wn\beta(\tr x)).$$
Here the left hand side is equivalent by ($\wn_\exists$) and \cite{M1}*{\ref{g1:move-oc-wn}(a)}
to $\oc\wn\alpha\to\oc\wn\exists\tr x\,\beta(\tr x)$, which in the presence of
$\prin\nabla\alpha\tofrom\alpha$ simplifies to $\alpha\to\exists\tr x\,\beta(\tr x)$.
On the other hand, $\exists\tr x\,(\alpha\to\beta(\tr x))$ can be rewritten in the presence of
$\prin\nabla\alpha\tofrom\alpha$ as $\exists\tr x\,(\oc\wn\alpha\to\oc\wn\beta(\tr x))$,
which by \cite{M1}*{\ref{g1:move-oc-wn}(a)} is equivalent to
$\exists\tr x\,\oc(\wn\alpha\to\wn\beta(\tr x))$.
Collecting the previous steps, we get that the Principle of Independence of Premise,
$\prin(\alpha\to\exists\tr x\,\beta(\tr x))\To \exists\tr x\,(\alpha\to\beta(\tr x))$, is equivalent 
(in the presence of the PC-principle) to
$\prin\oc\exists\tr x\,(\wn\alpha\to\wn\beta(\tr x))\To\exists\tr x\,\oc(\wn\alpha\to\wn\beta(\tr x))$;
and the assertion follows.
\end{proof}

\begin{proposition}\label{coupled principles}
(a) The principle $\prin\oc\exists\tr x\, p(\tr x)\Tofrom\exists\tr x\,\oc p(\tr x)$
\begin{enumerate}[label=\oldstylenums{\arabic*}.]
\item is equivalent to the conjunction of the $\exists_\Box$- and $\exists_\nabla$-principles; and
\item implies the Generalized Markov Principle.
\end{enumerate}

(b) The principle $\prin\wn\forall\tr x\,\alpha(\tr x)\Tofrom\forall\tr x\,\wn\alpha(\tr x)$
\begin{enumerate}[label=\oldstylenums{\arabic*}.]
\item is equivalent to the conjunction of the $\forall_\Box$- and $\forall_\nabla$-principles; and
\item implies the Parametric Distributivity Principle.
\end{enumerate}
\end{proposition}

The Parametric Distributivity Principle,
$$\prin\neg\forall\tr y\,(\alpha(\tr y)\lor\forall\tr x\,\beta(\tr x,\tr y))\To
\neg\forall\tr y\,\forall\tr x\,(\alpha(\tr y)\lor\beta(\tr x,\tr y)),$$
introduced in \cite{M0}*{Example \ref{int:constant domain}}, is a variation of a principle of Kleene; each
of the two principles is a double negation of a two-variable generalization of the Constant Domain Principle
(whereas the double negation of the one-variable Constant Domain Principle is intuitionistically derivable).

\begin{proof}[Proof. \o1's] We will check (b), \o1; an entirely similar argument proves (a), \o1.
We have, in bare QHC calculus,
$\turnstile\wn\forall\tr x\,\alpha(\tr x)\To\wn\forall\tr x\,\nabla\alpha(\tr x)$
and
$\turnstile\wn\forall\tr x\,\nabla\alpha(\tr x)\To\forall\tr x\,\wn\alpha(\tr x)$,
and by \ref{principles-q2}(a,c), the reverse implications are precisely
the $\forall_\Box$- and $\forall_\nabla$-principles.
Thus the meta-conjunction of the latter principles is equivalent to
$\prin\forall\tr x\,\wn\alpha(\tr x)\To\wn\forall\tr x\,\alpha(\tr x)$.
But the converse to the latter is nothing but the law ($\wn_\forall$).
\end{proof}

\begin{proof}[Proof. (a), \o2] By \cite{M1}*{\ref{g1:insolubility}} and ($\wn_\forall$) we have in QHC
$\turnstile\neg\forall\tr x\,\alpha(\tr x)\Tofrom\oc\neg\wn\forall\tr x\,\alpha(\tr x)$
and $\turnstile\oc\neg\wn\forall\tr x\,\alpha(\tr x)\To\oc\neg\forall\tr x\,\wn\alpha(\tr x)$.
Since the Generalized Markov Principle holds classically, we also have
$\turnstile\oc\neg\forall\tr x\,\wn\alpha(\tr x)\To\oc\exists\tr x\,\neg\wn\alpha(\tr x)$.
From the hypothesis we also get 
$\prin\oc\exists\tr x\,\neg\wn\alpha(\tr x)\To\exists\tr x\,\oc\neg\wn\alpha(\tr x)$.
Finally, $\turnstile \exists\tr x\,\oc\neg\wn\alpha(\tr x)\Tofrom\exists\tr x\,\neg\alpha(\tr x)$ by
\cite{M1}*{\ref{g1:insolubility}}.
Thus, using the hypothesis, we have obtained $\prin\neg\forall\tr x\alpha(\tr x)\to\exists\tr x\neg\alpha(\tr x)$,
the Generalized Markov Principle.
\end{proof}

\begin{proof}[(b), \o2] Since the Parametric Distributivity Principle holds classically, we have
$$\turnstile\oc\neg\forall\tr y\,(\wn\alpha(\tr y)\lor\forall\tr x\,\wn\beta(\tr x,\tr y))\To
\oc\neg\forall\tr y\,\forall\tr x\,(\wn\alpha(\tr y)\lor\wn\beta(\tr x,\tr y)).$$
Since $\wn$ commutes with $\lor$ and, by the hypothesis, with $\forall$, we get
$$\prin\oc\neg\wn\forall\tr y\,(\alpha(\tr y)\lor\forall\tr x\,\beta(\tr x,\tr y))\To
\oc\neg\wn\forall\tr y\,\forall\tr x\,(\alpha(\tr y)\lor\beta(\tr x,\tr y)).$$
Finally, by \cite{M1}*{\ref{g1:insolubility}} the latter is equivalent to the Parametric
Distributivity Principle.
\end{proof}

\begin{proposition}\label{Jankov}\

(a) The following are equivalent:
\begin{enumerate}[label=\oldstylenums{\arabic*}.]
\item $\prin\oc(p\lor q)\Tofrom\oc p\lor\oc q$;
\item $\prin\Dec(p)$;
\item the meta-conjunction of the K-principle, $\prin\Box(p\lor q)\tofrom\Box p\lor\Box q$, and the D-principle,
$\prin\nabla(\alpha\lor\beta)\tofrom\nabla\alpha\lor\nabla\beta$;
\item the meta-conjunction of the K-principle and Jankov's principle, $\prin\neg\alpha\lor\neg\neg\alpha$.
\end{enumerate}

(b) The following are equivalent:
\begin{enumerate}[label=\oldstylenums{\arabic*}.]
\item $\prin\Dec(\nabla\alpha)$;
\item the meta-conjunction of the H-principle and the D-principle;
\item the meta-conjunction of the H-principle and Jankov's principle.
\end{enumerate}
\end{proposition}

As noted above (see \ref{Hilbert}), $\prin\Dec(\nabla\alpha)$
is equivalent to $\prin\Dec(\Box p)$.

Beware that while $\prin\Dec(\nabla\alpha)$ is a special case of $\prin\alpha\lor\neg\alpha$, the principle
$\prin\Dec(p)$ does not follow from it (using the model of \ref{strict implications} below).

\begin{proof}[Proof. (b)] By \cite{M1}*{\ref{g1:move-nabla}}, \o1 is equivalent to
$\prin\nabla\alpha\lor\neg\alpha$.
By \cite{M1}*{\ref{g1:nabla-negneg1}}, $\turnstile\nabla\alpha\to\neg\neg\alpha$.
Hence \o1 implies $\prin\neg\neg\alpha\lor\neg\alpha$.

On the other hand, $\prin\Dec(\nabla\alpha)$ implies $\prin\wn\Dec(\nabla\alpha)$,
which by \ref{Hilbert} is equivalent to the H-principle.
Thus \o1 implies \o3.

Conversely, assuming the H-principle, we have $\prin\nabla\alpha\tofrom\neg\neg\alpha$, and our derivation of
$\prin\neg\neg\alpha\lor\neg\alpha$ from \o1 can be reversed.

Finally, assuming the H-principle, we also have that
$\prin\nabla(\alpha\lor\beta)\Tofrom\nabla\alpha\lor\nabla\beta$ is equivalent to
$\prin\neg\neg(\alpha\lor\beta)\Tofrom\neg\neg\alpha\lor\neg\neg\beta$.
But the latter is equivalent to $\prin\neg\alpha\lor\neg\neg\alpha$ (see \cite{M0}*{\ref{int:Jankov's logic}}).
\end{proof}

\begin{proof}[(a). \o1$\iff$\o2.] If we specialize $\prin\oc(p\lor q)\to\oc p\lor\oc q$ by substituting 
$\neg p$ for $q$, we obtain $\prin \oc p\lor\oc\neg p$ using the classical law of excluded middle.

Conversely, since $\oc$ commutes with $\land$, and $(p\lor q)\land\neg p$ is intuitionistically
equivalent to $q\land\neg p$, we have $\turnstile \oc(p\lor q)\land\oc\neg p\tofrom \oc q\land\oc\neg p$.
In particular, we obtain $\turnstile\oc(p\lor q)\land\oc\neg p\to\oc q$.
By the exponential law, this rewrites as $\turnstile\oc(p\lor q)\to(\oc\neg p\to\oc q)$.
This is easily seen to yield $\turnstile\oc(p\lor q)\to(\oc p\lor\oc\neg p\to\oc p\lor\oc q)$.
A double application of the exponential law transforms this into
$\turnstile\oc p\lor\oc\neg p\to\big(\oc(p\lor q)\to\oc p\lor\oc q\big)$.
\end{proof}

\begin{proof}[\o1$\iff$\o3] This is proved similarly to \ref{coupled principles}, using
\ref{principles}(b) and \ref{principles2}(b).
\end{proof}

\begin{proof}[\o3$\iff$\o4] This is similar to (b), using
that the K-principle implies the H-principle.
\end{proof}

\metameta

\subsection{ED-principle}
The following {\it Exclusive Disambiguation Principle} (ED-principle): 
\[\fm{\prin\neg(\alpha\land\beta)\To\big(\nabla(\alpha\lor\beta)\to\nabla\alpha\lor\nabla\beta\big)}\]
is motivated by the observations that (i)
$\nabla\Gamma$ has a solution if and only if $\Gamma$ does (in fact, 
$\fm{\turnstile \wn\nabla\gamma\tofrom\wn\gamma}$ in QHC),
and that (ii) $\nabla\Gamma=\oc(\wn\Gamma)$ can have at most one solution according to Lafont's
argument (see \cite{M0}*{\S\ref{int:Lafont}}).
Indeed, by (i), $\nabla\Gamma\lor\nabla\Delta$ has a solution if and only if $\nabla(\Gamma\lor\Delta)$
does; whereas by (ii), $\nabla\Gamma\lor\nabla\Delta$ can have at most one solution as long as
$\nabla\Gamma\land\nabla\Delta$ is known to have no solutions.
By (i), the latter assumption amounts to $\Gamma\land\Delta$ having no solutions.
Now under some circumstances, knowing existence {\it and} uniqueness of a solution could amount to
the knowledge of that solution itself.

\formulas
Clearly, the D-principle, $\prin\nabla(\alpha\lor\beta)\Tofrom(\nabla\alpha\lor\nabla\beta)$, implies the
ED-principle.

\begin{proposition}\label{ED-prin} The following are equivalent forms of the ED-principle:

(i) $\fm{\dfrac{\neg(\alpha\land\beta)}{\nabla(\alpha\lor\beta)\to\nabla\alpha\lor\nabla\beta}}$;
\medskip

(ii) $\prin\nabla\big((\alpha\land\neg\beta)\lor(\beta\land\neg\alpha)\big)\To\nabla\alpha\lor\nabla\beta$.
\end{proposition}

\begin{proof}
By the exponential law, the ED-principle is equivalent to
$\prin\neg(\alpha\land\beta)\land\nabla(\alpha\lor\beta)\To\nabla\alpha\lor\nabla\beta$.
To see that the latter is equivalent to (ii), let us note that by \cite{M1}*{\ref{g1:move-nabla} and
\ref{g1:move-box-diamond}}, $\neg(\alpha\land\beta)\land\nabla(\alpha\lor\beta)$ is equivalent to
$\nabla\big((\alpha\lor\beta)\land\neg(\alpha\land\beta)\big)$.
By \cite{M0}*{(\ref{int:morgan1}) and (\ref{int:implication-exp})},
the latter is in turn equivalent to
$\nabla\Big(\big(\alpha\land(\alpha\to\neg\beta)\big)\lor\big(\beta\land(\beta\to\neg\alpha)\big)\Big)$,
and the assertion follows.

Now, the ED-principle implies (i) using the modus ponens.
To see that (i) implies (ii), we substitute $\Phi:=\gamma\land\neg\delta$ for $\alpha$ and
$\Psi:=\delta\land\neg\gamma$ for $\beta$, and observe that 
$\turnstile\neg\big((\gamma\land\neg\delta)\land(\gamma\land\neg\delta)\big)$.
Hence (i) implies $\prin\nabla(\Phi\lor\Psi)\to\nabla\Phi\lor\nabla\Psi$, which in turn implies
(ii) using that
$\turnstile\nabla(\gamma\land\neg\delta)\lor\nabla(\delta\land\neg\gamma)\to\nabla\gamma\lor\nabla\delta$.
\end{proof}

\begin{proposition}\label{Jankov2}
(a) Jankov's principle, $\prin\neg\alpha\lor\neg\neg\alpha$, implies the ED-principle.

(b) $\prin\Dec(\nabla\alpha)$ is equivalent to the meta-conjunction of the H- and ED-principles.

(c)  $\prin\Dec(p)$ is equivalent to the meta-conjunction of the K- and ED-principles.
\end{proposition}

\begin{proof}[Proof. (a)] Let us recall that $\prin\neg\gamma\lor\neg\neg\gamma$ implies de Morgan's
principle $\prin\neg(\gamma\land\delta)\to\neg\gamma\lor\neg\delta$ (the converse of an
intuitionistically derivable principle), see \cite{M0}*{\ref{int:Jankov's logic}}.
Thus it suffices to show that the formula
$\neg\alpha\lor\neg\beta\to\big(\nabla(\alpha\lor\beta)\to\nabla\alpha\lor\nabla\beta\big)$ is derivable in QHC.
By the exponential law, it is equivalent to 
$\nabla(\alpha\lor\beta)\to(\neg\alpha\lor\neg\beta\to\nabla\alpha\lor\nabla\beta)$.
By \cite{M0}*{(\ref{int:dm1})}, the latter is equivalent to
$\nabla(\alpha\lor\beta)\to(\neg\alpha\to\nabla\alpha\lor\nabla\beta)
\land(\neg\beta\to\nabla\alpha\lor\nabla\beta)$,
which by \cite{M0}*{(\ref{int:anti-dm1})} is in turn equivalent to $\mm{F\land G}$, where
$\mm F=\nabla(\alpha\lor\beta)\to(\neg\alpha\to\nabla\alpha\lor\nabla\beta)$ and 
$\mm G=\nabla(\alpha\lor\beta)\to(\neg\beta\to\nabla\alpha\lor\nabla\beta)$.
Here $\mm F$ is equivalent, using the exponential law and \cite{M1}*{\ref{g1:move-nabla}}, to
$\nabla(\alpha\lor\beta)\land\nabla\neg\alpha\to\nabla\alpha\lor\nabla\beta$.
The latter is equivalent, using \cite{M1}*{\ref{g1:move-box-diamond}} and
\cite{M0}*{(\ref{int:morgan1})}, to
$\nabla\big((\alpha\land\neg\alpha)\lor(\beta\land\neg\alpha)\big)\to\nabla\alpha\lor\nabla\beta$.
This is derivable in QHC, since $\turnstile\nabla(\beta\land\neg\alpha)\to\nabla\beta$.
Similarly, $\mm G$ is also derivable in QHC, and the assertion follows.
\end{proof}

\begin{proof}[(b)] One implication follows from (a) and \ref{Jankov}(b).
Conversely, the H-principle implies, in particular, $\prin\wn\Dec(\nabla\gamma)$, that is,
$\prin\nabla(\nabla\gamma\lor\neg\nabla\gamma)$.
Substituting $\nabla\gamma$ and $\neg\nabla\gamma$ for $\alpha$ and $\beta$ in the ED-principle and
using \cite{M1}*{\ref{g1:move-nabla}}, we get $\prin\nabla\gamma\lor\neg\nabla\gamma$,
i.e., $\prin\Dec(\nabla\gamma)$.
\end{proof}

\begin{proof}[(c)] One implication follows from \ref{Jankov}(a).
The converse is parallel to (b).
In more detail, by \ref{Kolmogorov}, the K-principle is equivalent to $\prin\wn\Dec(p)$, that is,
$\prin\nabla(\oc p\lor\oc\neg p)$.
We have $\turnstile\oc p\land\oc\neg p\tofrom\oc(p\land\neg p)$
and $\turnstile\oc(p\land\neg p)\to\oc \clbot$ and $\turnstile\oc \clbot\tofrom\ab$.
Hence, by substituting $\oc p$ and $\oc\neg p$ for $\alpha$ and $\beta$ in the ED-principle, we get
$\prin\oc p\lor\oc\neg p$, i.e., $\prin\Dec(p)$.
\end{proof}

\metameta

\section{Topological models}\label{models}

\subsection{Interior-based models} \label{interior-based}
By composing the $\Box$-interpretation of QHC (see \cite{M1}) with a topological model of QS4 
(see \cite{M0}*{\S\ref{int:S4}}), we obtain an {\it interior-based model} of QHC
(also known as an ``Euler--Tarski model'').
In relatively self-contained terms, such a model can be described as follows:
\begin{itemize}
\item We fix a topological space $X$ and a domain of discourse $\D$.
\item The intuitionistic side is interpreted as in the Tarski model (see \cite{M0}*{\S\ref{int:Tarski}})
in open subsets of $X$.
\item The classical side is interpreted as in the Leibniz--Euler model (see \cite{M0}*{\S\ref{int:Euler}})
in arbitrary subsets of $X$.
\item $|\wn\Phi|=|\Phi|$, and $|\oc F|=\Int|F|$.
\end{itemize}
Even though $\wn$ is interpreted by the identity, it need not commute with $\to$ or with $\forall$
in an interior-based model, since the interpretations of intuitionistic $\to$ and $\forall$ do not coincide 
with the restrictions of the interpretations of classical $\to$ and $\forall$.
Indeed, the interpretations of intuitionistic $\to$ and $\forall$ can be expressed in terms of
those of classical $\to$ and $\forall$ due to the validity of the following principles in all 
interior-based models:
\formulas
\begin{itemize}
\item[] $\prin(\alpha\to \beta)\Tofrom \oc(\wn\alpha\to\wn\beta)$;
\item[] $\prin\forall\tr x\,\alpha(\tr x)\Tofrom \oc\forall\tr x\, \wn\alpha(\tr x)$.
\end{itemize}
It should be noted that each of these two principles is equivalent to the PC-principle,
$\prin\nabla\alpha\tofrom\alpha$ (see \ref{expressive}(a) and \ref{principles-q}(a)).
\metameta

Now that we have a model of the QHC calculus, we obtain

\begin{theorem} The QHC calculus is consistent.
\end{theorem}

It is easy to see that in an interior-based model, decidable propositions
are represented by clopen sets, and stable propositions by sets $S$
such that $\Int\Cl S=\Int S$ (these include closed sets as well as
regular open sets).

\subsection{Regularization-based models} \label{regularization-based}
By composing the $\neg\neg$-interpretation of QHC (see \cite{M1}) with a Tarski model of QH 
(see \cite{M0}*{\S\ref{int:Tarski}}), we get a {\it regularization-based model} of QHC
(also known as a ``Tarski--Kolmogorov model'').
A relatively self-contained description is as follows:
\begin{itemize}
\item We fix a topological space $X$ and a domain of discourse $\D$.
\item The intuitionistic side is interpreted as in a Tarski model in open sets of $X$.
\item The classical side is interpreted in regular open sets of $X$, via the $\neg\neg$-translation.
\item $|\oc F|=|F|$, and $|\wn\Phi|=\Int(\Cl|\Phi|)$.
\end{itemize}
In other words,  $|\wn\Phi|=|\neg\neg\Phi|$.
On the classical side we have, in more detail:
\begin{itemize}
\item $|p|(x_1,\dots,x_n)$ is a regular open set for each $n$-ary predicate variable $p$;
\item $|F\land G|=|F|\cap|G|$;
\item $|\forall\tr x\, H(\tr x)|=\Int\bigcap_x |H(\tr x)|$;
\item$|\neg F|=\Int(X\but |F|)$;
\item other classical connectives reduce to $\land$, $\forall$ and $\neg$; namely,
\item $|F\lor G|=\Int\Cl(|F|\cup|G|)$;
\item $|\exists\tr x\, H(\tr x)|=\Int\Cl\bigcup_x |H(\tr x)|$;
\item $|F\to G|=\Int\big((X\but |F|)\cup|G|\big)$;
\item $|\top|=X$ and $|\bot|=\emptyset$,
\end{itemize}
where $F$ and $G$ are c-formulas and $H$ is a 1-c-formula, and $\iass$ and $\pval$ are fixed.

Even though $\oc$ is interpreted by the identity, it need not commute
with $\lor$ or with $\exists$ in a regularization-based model, since
the interpretations of classical $\lor$ and $\exists$ do not coincide with the restrictions of the
interpretations of intuitionistic $\lor$ and $\exists$.
Indeed, the interpretations of classical $\lor$ and $\exists$ are expressible in terms of
those of intuitionistic $\lor$ and $\exists$ due to the validity of the following principles in all 
regularization-based models:
\formulas
\begin{itemize}
\item[] $\prin p\lor q\Tofrom\wn(\oc p\lor \oc q)$;
\item[] $\prin\exists\tr x\, p(\tr x)\Tofrom \wn\exists\tr x\, \oc p(\tr x)$.
\end{itemize}
It should be noted that each of these two principles is equivalent to the K-principle,
$\prin\Box p\tofrom p$ (see \ref{expressive}(b) and \ref{principles-q}(b)).
\metameta

Let us note that $\Diamond$ corresponds to the closure operator in interior-based models of QHC (or 
in topological models of QS4), and the subset $\Int\Cl(\Int S)$ of a topological space is precisely 
the smallest regular open set containing the open set $\Int S$.
It follows from this that the $\Diamond$-interpretation  $A\mapsto A_\Diamond$ 
(see \cite{M1}*{\S\ref{g1:diamond}}) ``pulls back'' a regularization-based model out of 
every interior-based model:

\begin{proposition} Every interior-based model $M$ determines a regularization-based model $M_\Diamond$ 
as follows: the interpretation of a formula $A$ of QHC in $M_\Diamond$ is $|A|_{M_\Diamond}:=|A_\Diamond|_M$.
\end{proposition}

\begin{remark} \label{Diamond-meta}
By \cite{M1}*{Theorem \ref{g1:Diamond-interpretation}}, $\turnstile A$ implies
$\turnstile A_\Diamond$ if $A$ is a formula of QHC.
This meta-judgement cannot be internalized.
Indeed, $A\turnstile A_\Diamond$ would imply that principles satisfied in interior-based models
are also satisfied in regularization-based models --- which is not the case, as we will see in
Theorem \ref{intermediate models}.
\end{remark}

\subsection{Examination of principles}

\begin{theorem} \label{intermediate models}
(a) The R-principle is valid in all interior-based models and in all regularization-based models.

(b) The PC-, D- and ED-principles; the PC-rule; and the $\forall_\nabla$- and
$\exists_\nabla$-principles are valid in all interior-based models but not valid in some 
regularization-based models.

(c) The K- and H-principles and the $\forall_\Box$- and $\exists_\Box$-principles
are valid in all regularization-based models but not valid in some interior-based models.
\end{theorem}

\begin{proof} By the results of \S\ref{principles section}, the principles listed in (a) and (b) all 
follow from $\fm{\prin\alpha\tofrom\nabla\alpha}$, which is clearly valid in all interior-based models;
whereas the principles listed in (a) and (c) all follow from $\fm{p\tofrom\Box p}$,
which is clearly valid in all regularization-based models.

Thus it remains to check the negative assertions, of which those on the D- and K-principles
can be omitted, since they imply respectively the ED- and H-principles.

In the presence of the K-principle, which is valid in all regularization-based models,
the PC-principle is equivalent to the PC-rule, which is in turn equivalent to
$\fm{\prin\alpha\lor\neg\alpha}$; the D- and ED-principles are equivalent to
$\fm{\prin\neg\alpha\lor\neg\neg\alpha}$; and the $\forall_\nabla$- and $\exists_\nabla$-principles
are equivalent respectively to the $\neg\neg$-Shift Principle and the Strong Markov Principle.
But each of these is not valid in some Tarski model of QH (see \cite{M0}*{\S\ref{int:topology}}),
which in turn extends to a regularization-based model of QHC.

Similarly, in the presence of the PC-principle, which is valid in all interior-based models,
the H-principle is equivalent to $\fm{\prin\alpha\lor\neg\alpha}$; and the $\forall_\Box$- and the
$\exists_\Box$-principles imply respectively the Constant Domain Principle and the
Principle of Independence of Premise.
But each of these is not valid in some Tarski model of QH (see \cite{M0}*{\S\ref{int:topology}}),
which in turn extends to an interior-based model of QHC.
\end{proof}

\begin{proposition} \label{strict implications}
(a) The H-principle implies neither the K-principle nor the $\exists_\Box$-principle.

(b) The D-principle, $\fm{\prin\nabla(\alpha\lor\beta)\tofrom\nabla\alpha\lor\nabla\beta}$,
does not imply the PC-rule, $\fm{\dfrac{\nabla\alpha}\alpha}$.
\end{proposition}

\begin{proof}[Proof. (a)]
Let us consider an interior-based model on $\R$ with respect to the valuation field where all problem variables 
are interpreted by $\emptyset$ or $\R$ on every input, and all predicate variables are interpreted by $\emptyset$ 
or $\Q$ or $\R\but\Q$ or $\R$ on every input.
(This is possible since both $\Q$ and $\R\but\Q$ have empty interior.)
Since all problems are interpreted trivially, the H-Principle (e.g.\ in the form
$\fm{\prin\nabla\alpha\tofrom\neg\neg\alpha}$) is valid.
But the interior operator is nontrivial, so the K-Principle (e.g.\ in the form
$\fm{\prin\Box p\tofrom p}$) is not valid.
The $\exists_\Box$-principle also fails, e.g.\ by considering a two-element domain.
\end{proof}

\begin{proof}[(b)] Let us consider the regularization-based model on the two-element poset $0<1$,
viewed as an Alexandroff space.
The regularization of the open set $\{1\}$ is the entire space, $\{0,1\}$, so
$\fm{\nabla\alpha\not\Turnstile\alpha}$.
On the other hand, using that the only regular open sets are $\emptyset$ and $\{0,1\}$, it is easy to
check that $\fm{\Turnstile\nabla(\alpha\lor\beta)\tofrom\nabla\alpha\lor\nabla\beta}$.
\end{proof}

\begin{remark} The PC-rule implies the PC-principle (and in particular the weaker D-, $\forall_\nabla$- 
and $\exists_\nabla$-principles) both in interior-based models and in regularization-based models, because 
all of them satisfy the R-principle (see \ref{R+PC}).
\end{remark}

\subsection{Subset/sheaf-valued models} \label{sheaf model}

A {\it subset/sheaf-valued structure} interpreting the language of QHC is described as follows.

\begin{enumerate}[label=(\roman*)]
\item We fix a topological space $B$ and a domain of discourse $\D$.
\item The intuitionistic side is interpreted as in a sheaf-valued model of \cite{M0}*{\S\ref{int:sheaves}} 
over $B$.
Thus, nullary problem variables are interpreted by sheaves (of sets) over $B$; $n$ary problem variables
by families of sheaves indexed by $\D^n$; and intuitionistic connectives and quantifiers by 
the usual operations on sheaves:
\begin{itemize}
\item $|\Phi\lor\Psi|=|\Phi|\sqcup|\Psi|$;
\item $|\Phi\land\Psi|=|\Phi|\x|\Psi|$;
\item $|\Phi\to\Psi|=\Hom(|\Phi|,|\Psi|)$;
\item $|\ab|=\Char\emptyset$;
\item $|\exists \tr x\,\Xi(\tr x)|=\bigsqcup_{d\in\D}|\Xi|(d)$;
\item $|\forall \tr x\,\Xi(\tr x)|=\prod_{d\in\D}|\Xi|(d)$,
\end{itemize}
where $\Phi$ and $\Psi$ are i-formulas and $\Xi$ is a 1-i-formula, and $\iass$ and $\pval$ are fixed.
Also, if $\Phi$ is an $n$-i-formula, $\Turnstile\Phi$ means that each sheaf in the family of sheaves 
$|\Phi|$ has a global section.

\item The classical side is interpreted as in a Leibniz--Euler model (see \cite{M0}*{\S\ref{int:Euler}}) 
--- in arbitrary subsets of $B$. Thus:
\begin{itemize}
\item $|F\lor G|=|F|\cup|G|$;
\item $|F\land G|=|F|\cap|G|$;
\item $|\neg F|=X\but |F|$;
\item $|F\to G|=(X\but |F|)\cup |G|$;
\item $|\exists \tr x\,H(\tr x)|=\bigcup_{d\in\D}|H|(d)$;
\item $|\forall \tr x\,H(\tr x)|=\bigcap_{d\in\D}|H|(d)$,
\end{itemize}
where $F$ and $G$ are c-formulas and $H$ is a 1-c-formula, and $\iass$ and $\pval$ are fixed.
Also, if $F$ is an $n$-c-formula, $\Turnstile F$ means that each subset in the family of 
subsets $|F|$ coincides with $X$.

\item $|\oc F|=\chi_{\Int|F|}$, the sheaf provided by the inclusion $\Int|F|\emb B$.
\item $|\wn\Phi|=\Supp|\Phi|$, the set of all $b\in B$ such that the stalk $|\Phi|_b\ne\emptyset$.
\end{enumerate}

In exactly the same way we define subset/{\it pre\/}sheaf-valued structures.
The only essential difference between subset/sheaf- and subset/presheaf-valued structures is in 
the interpretations of intuitionistic $\lor$ and $\exists$, due to the fact that 
$\sigma\F\sqcup\sigma\G\ne\sigma(\F\sqcup\G)$ in general.
We write $\Char\F$ for the presheaf of sections of the characteristic sheaf $\chi_U$, provided by
the inclusion of the open subset $U\emb B$.

As a simple illustration, let us discuss the difference between the interpretations of $\nabla$ and 
$\neg\neg$ in a subset/presheaf-valued model.
Since $\Supp(F)$ is always open for a presheaf $F$, $|\oc\wn\Phi|=\Char V$, where $V=\Supp|\Phi|$.
On the other hand, by \cite{M0}*{\ref{int:sheaf negation}}, $|\neg\Phi|=\Char U$, where $U=\Int(B\but V)$.
Hence $|\neg\neg\Phi|=\Char V'$, where $V'=\Int(B\but U)=\Int(\Cl V)$.

\begin{theorem} \label{global1} Subset/sheaf- and subset/presheaf-valued structures are models of QHC.
\end{theorem}

\begin{proof} It is well-known that subsets interpret classical predicate logic (see \cite{RS}).
It is also known that sheaves and presheaves interpret intuitionistic
predicate logic (see \cite{M0}*{\S\ref{int:sheaves}}).
Thus it remains to check that the additional laws and inference rules of QHC are
satisfied in subset/(pre)sheaf-valued structures.

\begin{enumerate}
\item[($\wn_\top$)] $\dfrac{\fm\alpha}{\wn\fm\alpha}$;
\medskip
\item[($\oc_\top$)] $\dfrac{\fm p}{\oc\fm p}$.
\end{enumerate}
\smallskip

Indeed, if $|\fm\alpha|$ has a global section, then $\Supp|\fm\alpha|=B$.
Conversely, if $|\fm p|=B$, then $\Char(\Int|\fm p|)$ has a global section.

\begin{enumerate}[resume,label=(\roman*)]
\item[($\wnto$)] $\prin\fm{\wn(\alpha\to \beta)\To(\wn\alpha\to\wn\beta)}$
\end{enumerate}

This holds since $\Supp\Hom(F,G)\subset\Int\big(\Supp G\cup (B\but\Supp F)\big)$ (see
\cite{M0}*{\ref{int:BHK-verification}}).

\smallskip

\begin{enumerate}
\item[($\octo$)] $\prin\fm{\oc(p\to q)\To(\oc p\to \oc q)}$.
\end{enumerate}

We have $|\fm{\oc(p\to q)}|=\Char U$, where $U=\Int|\fm{p\to q}|=\Int\big((B\but|\fm p|)\cup|\fm q|\big)$.
On the other hand, $|\fm{\oc p\to\oc q}|=\Hom\big(\Char(\Int|\fm p|),\Char(\Int|\fm q|)\big)$.
By \cite{M0}*{\ref{int:Hom-sheaf}(a)}, the latter presheaf is isomorphic
to $\Char V$, where $V=\Int\big((B\but\Int|\fm p|)\cup\Int|\fm q|\big)$.

To construct a presheaf morphism $\Char U\to\Char V$, it suffices to show that $U\subset V$.
If $b\in U$, then $b$ has an open neighborhood $W$ in $B$ such that $W\subset (B\but|\fm p|)\cup|\fm q|$.
Then $W\cap|\fm p|\subset W\cap|\fm q|$.
Hence $\Int_W(W\cap|\fm p|)\subset\Int_W(W\cap|\fm q|)$.
Since $W$ is open, $\Int_W(W\cap S)=W\cap\Int S$ for any $S\subset B$.
Then $W\cap\Int|\fm p|\subset W\cap\Int|\fm q|$.
Hence $W\subset (B\but\Int|\fm p|)\cup\Int|\fm q|$.
Thus $b\in V$.

\begin{enumerate}
\item[($\oc_\bot$)] $\neg\oc \clbot$.
\end{enumerate}

Since $\turnstile \clbot\tofrom\fm p\land\neg\fm p$, we have $|\clbot|=\emptyset$, so 
$|\oc \clbot|$ is the empty sheaf.
Hence $|\oc \clbot\to\ab|=\Hom(\Char\emptyset,\Char\emptyset)=\Char B$.

\begin{enumerate}
\item[($\wn\oc$)] $\prin\wn\oc\fm p\to\fm p$.
\end{enumerate}

We have $|\wn\oc\fm p|=\Int|\fm p|$, so $|\wn\oc\fm p\to\fm p|=|\neg\wn\oc\fm p\lor\fm p|=
(B\but\Int|\fm p|)\cup|\fm p|=B$.

\begin{enumerate}
\item[($\oc\wn$)] $\prin\fm\alpha\to \oc\wn\fm\alpha$.
\end{enumerate}

Since $\Supp(F)$ is always open, $|\oc\wn\fm\alpha|=\Char V$, where $V=\Supp|\fm\alpha|$.
In other words, $V$ is the union of all open sets $U$ such that $|\fm\alpha|(U)\ne\emptyset$.
Thus $|\oc\wn\fm\alpha|(U)=\{\id_U\}$ if $|\fm\alpha|(U)\ne\emptyset$, and else $|\oc\wn\fm\alpha|(U)=\emptyset$.
Then the constant maps $\phi(U)\:|\fm\alpha|(U)\to|\oc\wn\fm\alpha|(U)$, $s\mapsto\id_U$, combine into
a natural transformation $\phi\:|\fm\alpha|\to|\oc\wn\fm\alpha|$.
\end{proof}

Interior-based models of QHC are not the same as subset/(pre)sheaf-valued models of QHC where all problem 
variables are interpreted by the characteristic (pre)sheaves of open subsets of $B$.
Indeed, formulas containing $\lor$ and $\exists$ will generally not be interpreted by such (pre)sheaves, 
but only by their coproducts.
In fact, the connection between subset/(pre)sheaf and interior-based models is provided
by the $\nabla$-translation, $A\mapsto A_\nabla$ (see \cite{M1}*{\S\ref{g1:Euler-Tarski1}}):

\begin{proposition} Every subset/(pre)sheaf-valued model $M$ determines an interior-based model $M_\nabla$: 
the interpretation of a formula $A$ of QHC in $M_\nabla$ is $|A|_{M_\nabla}:=|A_\nabla|_M$.
\end{proposition}

\begin{remark}\label{nabla-meta}
By \cite{M1}*{Theorem \ref{g1:nabla-interpretation}}, if $\Phi$ is a formula of intuitionistic logic, 
$\turnstile\Phi$ if and only if $\turnstile\Phi_\nabla$.
This meta-judgement cannot be internalized.
Indeed, for $\Phi=\forall\tr x\,\exists\tr y\, \fm\alpha(\tr x,\tr y)$ it follows from
\cite{M0}*{Proposition \ref{int:Tarski forall fails}} that $\Phi_\nabla\not\turnstile\Phi$ and also
$\neg\Phi\not\turnstile (\neg\Phi)_\nabla$.
\end{remark}

\begin{remark}
The class of models of QHC pulled back from sheaf models via the $\Diamond$-translation
can be described as follows in closed terms:

\begin{itemize}
\item We fix a topological space $B$ and a domain $\D$.
\item The intuitionistic side is interpreted as in the sheaf-valued model in $B$, with $\D$;
\item The classical side is interpreted as in a regularization-based model in $B$, with $\D$;
\item $|\oc F|=\chi_{|F|}$, the sheaf provided by the inclusion $|F|\emb B$.
\item $|\wn\Phi|=\Int\Cl(\Supp|\Phi|)$, the regularization of $\Supp|\Phi|$.
\end{itemize}

These models satisfy the K-principle (e.g.\ in the form $\prin\Box\fm p\tofrom\fm p$) and so do not seem to add
much new with respect to regularization-based models.
\end{remark}

\subsection{Principles revisited}

\begin{theorem}\label{sheaf-refutation}
(a) The R-principle (and hence also the H-, K- and PC-principles); the D-principle and
the PC-rule; and the $\forall_\Box$-, $\forall_\nabla$-, $\exists_\Box$- and
$\exists_\nabla$-principles are not valid in some subset/sheaf-valued models.

(b) The ED-principle is valid in all subset/sheaf-valued models.
\end{theorem}

\begin{proof}[Proof. (a)] The PC-rule is trivial to refute; for instance, the sheaf
$\chi_{(-\infty,1)}\sqcup\chi_{(-1,\infty)}$ over $\R$ has no global sections, yet its support is
the entire $\R$.

The D- and $\exists_\nabla$-principles assert that $\nabla$ commutes with
intuitionistic $\lor$ and $\exists$.
This is not the case for the characteristic sheaves of any non-disjoint open subsets of $B$ not
contained one in another.

The $\exists_\Box$- and $\forall_\Box$-principles are written essentially in the language of QS4 and
so can be refuted just like in the proof of Theorem \ref{intermediate models}.

Models where the $\forall_\nabla$- and R-principles fail are given by
\cite{M0}*{\ref{int:Tarski forall fails} and \ref{int:Tarski implication fails2}}, respectively.
\end{proof}

\formulas
\begin{proof}[(b)] Since $\turnstile\oc\neg\wn\Phi\tofrom\neg\Phi$, the equivalent rule form \ref{ED-prin}(i)
of the ED-principle can be obtained from $\dfrac{\oc\neg(p\land q)}{\oc(p\lor q)\to\oc p\lor\oc q}$ by 
substituting $p$, $q$ with $\wn\alpha$, $\wn\beta$.
The latter rule is equivalent to $\dfrac{\neg(p\land q)}{\oc(p\lor q)\to\oc p\lor\oc q}$.
We will show that $\neg(p\land q)\imp\oc(p\lor q)\to\oc p\lor\oc q$ is valid in any subset/sheaf-valued model
with respect to any valuation $\pval$ such that the sets $|p|^\pval$ and $|q|^\pval$ are open.
Indeed, $\neg(p\land q)$ is valid with respect to some $\pval$ if and only if the sets $|p|^\pval$ and 
$|q|^\pval$ are disjoint.
But given that they are disjoint, the sheaves $|\oc(p\lor q)|^\pval$ and $|\oc p|^\pval\sqcup|\oc q|^\pval$ 
are clearly isomorphic.
\end{proof}
\metameta

\begin{proposition} \label{some implic}
(a) The R-principle follows neither from the PC-rule nor from the D-principle.

(b) The PC-rule implies neither the D-principle nor the $\exists_\nabla$-principle.

(c) The PC-rule does not imply the $\forall_\nabla$-principle.
\end{proposition}

\begin{proof}[Proof. (a)] Let us consider the subset/sheaf-valued model of QHC in $\R^2$ with some $\D$ 
and let us fix the valuation $\pval$ that assigns to each $n$-ary problem variable and each $n$-tuple 
of elements of $\D$ the sheaf $\G$ of Example \cite{M0}*{\ref{int:Tarski implication fails2}} (the nontrivial 
double covering over $\R^2\but\{0\}$) and to each $n$-ary predicate variable and each $n$-tuple of elements of
$\D$ the subset $\emptyset$ of $\R^2$.
Then the R-principle is not valid with respect to the valuation field $\left<\pval\right>$, since for 
$F=\fm{\nabla\alpha}$ and $G=\fm\alpha$, 
$|\nabla F\to\nabla G|^\pval=\chi_{\R^2}$, whereas $|\nabla(F\to G)|^\pval$ has empty stalk at $0$.
On the other hand, for any formula $A$, its interpretation $|A|^\pval_\iass$ with respect to $\pval$ and 
some variable assignment $\iass$ is obtainable from $\G$ and $\emptyset$ by the operations that interpret 
the connectives of QHC (sheaf-theoretic operations, set-theoretic operations, and the operators 
$\F\mapsto\Supp\F$ and $S\mapsto\chi_S$).
Let us note that $\G$ and $\emptyset$ are invariant under linear transformations of $\R^2$.
Hence so is any sheaf of the form $|\Phi|^\pval_\iass$.
Consequently, for any valuation $\pval'\in\left<\pval\right>$, any assignment $\iass$ and any
problem variable $\phi$, $|\phi|^{\pval'}_\iass$ is invariant under linear transformations of $\R^2$.
Hence also any sheaf $\F$ of the form $|\Phi|^{\pval'}_\iass$ is invariant under linear transformations of $\R^2$.

If $\F_0\ne\emptyset$, then $\F$ has a section over some neighborhood of $0$ in $\R^2$.
But then by the scaling invariance $\F$ must have a global section.
Thus, in particular, $\fm{\nabla\alpha\,/\,\alpha}$ is valid with respect to $\left<\pval\right>$.

If $\F_p\ne\emptyset$ for some $p\ne 0$, then by the scaling and rotation invariance $\F_q\ne\emptyset$ 
for all $q\ne 0$.
Thus there are only three possibilities for $\Supp\F$: either $\emptyset$, or $\R^2$, or $\R^2\but\{0\}$.
Then for any other sheaf $\F'$ of the form $|\Phi|^{\pval'}_\iass$, either $\Supp\F\subset\Supp\F'$ or 
$\Supp\F'\subset\Supp\F$. 
Hence $\fm{\prin\nabla(\alpha\lor\beta)\to\nabla\alpha\lor\nabla\beta}$ is valid with respect to 
$\left<\pval\right>$.
\end{proof}

\begin{proof}[(b)] Let us consider the subset/sheaf-valued model of QHC in $\R^2$ with some $\D$ 
and let us fix a valuation $\pval$ that assigns to each $n$-ary predicate variable and each 
$n$-tuple of elements of $\D$ one of the subsets $U=\{(x,y)\mid x>0\}$ and $V=\{(x,y)\mid y>0\}$ of $\R^2$
and to each $n$-ary problem variable and each $n$-tuple of elements of $\D$ the sheaf $\chi_\emptyset$.
Since $U$, $V$ and $\chi_\emptyset$ are invariant under scaling of $\R^2$, so is the interpretation of 
any formula.
Similarly to the proof of (a), we get that $\fm{\nabla\alpha/\alpha}$ is valid with respect to $\left<\pval\right>$.
On the other hand, there exist no sheaf morphisms $\Char(U\cup V)\to\Char U\sqcup\Char V$, so
$\fm{\nabla(\alpha\lor\beta)\to\nabla\alpha\lor\nabla\beta}$ is not valid with respect to an appropriately chosen
$\pval$, and thus $\fm{\prin\nabla(\alpha\lor\beta)\to\nabla\alpha\lor\nabla\beta}$ is not valid with respect to 
$\left<\pval\right>$.
If we consider a two-element domain $\D$, we similarly get that the $\exists_\nabla$-principle is not valid
with respect to $\left<\pval\right>$, using an appropriately chosen $\pval$ of the above form.
\end{proof}

\begin{proof}[(c)] Let $U=\{re^{i\phi}\in\C\mid r>0,\,0<\phi<2\pi\}$ and
$V_n=\{re^{i\phi}\in\C\mid r>0,\,-\frac1n<\phi<\frac1n\}$.
Let us consider the subset/sheaf-valued model of QHC in $\R^2$ with $\D=\N$ 
and let us fix the valuation $\pval$ that assigns to each unary problem variable the family of sheaves
$\Char U\sqcup\Char V_n$, $n\in\N$, and to each $n$-ary predicate variable, $n\ne 1$, and each 
$n$-tuple of elements of $\D$ the sheaf $\chi_\emptyset$, and to each $n$-ary predicate variable and 
each $n$-tuple of elements of $\D$ the subset $\emptyset$ of $\R^2$.
Since $U$, $V_1,v_2,\dots$ and $\chi_\emptyset$ are invariant under scaling of $\C$ (i.e., multiplication
by a real number), so is the interpretation of any formula.
Hence it follows as in (a) that $\fm{\nabla\alpha/\alpha}$ is valid with respect to $\left<\pval\right>$.

On the other hand, for the unary problem variable $\fm\beta$ we have $|\nabla\fm\beta|^\pval(n)=\Char(\C\but\{0\})$ 
for each $n$, hence also $|\forall\tr x\,\nabla\fm\beta(\tr x)|^\pval=\Char(\C\but\{0\})$; whereas 
$|\nabla\forall\tr x\,\fm\beta(\tr x)|^\pval=\Char U$.
Thus the $\forall_\nabla$-principle is not valid with respect to $\left<\pval\right>$.
\end{proof}

\subsection{Summary}

\begin{theorem} (a) There are no implications between the K-, H-, R-, D-, PC-, $\forall_\Box$-, $\forall_\nabla$-,
$\exists_\Box$- and $\exists_\nabla$-principles and the PC-rule (denoted PC$^*$), 
except for the self-implications and the following implications:
\smallskip

\begin{center}
\includegraphics[width=8cm]{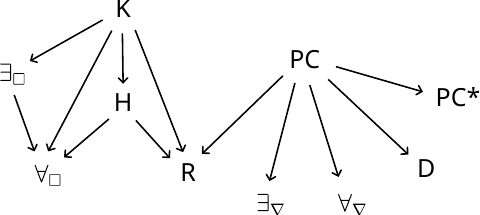}
\end{center}

(b) The ED-principle is not implied by any of the principles in (a), except the D-principle and the PC-principle.
(This does not say anything about the PC-rule.)
\end{theorem}

\begin{proof} The PC-principle, $\prin\nabla\alpha\tofrom\alpha$, clearly implies the PC-rule 
$\nabla\alpha\,/\,\alpha$, the R-principle $\prin(\nabla\alpha\to\nabla\beta)\to \nabla(\alpha\to\beta)$, 
the D-principle $\prin\nabla(\alpha\lor\beta)\to (\nabla\alpha\lor\nabla\beta)$, the $\forall_\nabla$-principle
$\prin\nabla\forall\tr x\, \alpha(\tr x)\tofrom \forall\tr x\, \nabla \alpha(\tr x)$ and the
$\exists_\nabla$-principle $\prin\exists\tr x\,\nabla\alpha(\tr x)\tofrom\nabla\exists\tr x\,\alpha(\tr x)$.
Also clearly the D-principle implies the ED-principle.

By \ref{Kolmogorov implies Hilbert} and \ref{Kolmogorov}, the K-principle implies the H-principle.
By \ref{HK-nabla-box} and \ref{h-quant} the H-principle implies the R-principle and the 
$\forall_\Box$-principle; and by \ref{quant-quant}, the $\exists_\Box$-principle implies 
the $\forall_\Box$-principle.

It remains to show that there are no other implications.
By Theorem \ref{intermediate models}, each of the principles/rules in (1) implies none of the principles in (2),
except possibly the R-principle, which is in (1); and each of the principles in (2) implies none of 
the principles/rules in (1), except possibly the R-principle, which is in (2).
Then it remains to show that there no extra implications between the principles/rules in (1), and no extra 
implications between the principles in (2).
Also, by Theorem \ref{intermediate models}, the R-principle implies no other principles/rules in (1) and (2).

By Theorem \ref{sheaf-refutation}, the ED-principle implies no other principles/rules in (1).
By \ref{no-implic}, \ref{strict implications}(b) and \ref{some implic}(a), the D-principle does not imply
any of the R-, $\exists_\nabla$- and $\forall_\nabla$-principles and the PC-rule.
Hence it also does not imply the PC-principle, and thus it implies only the principles in (4).
By \ref{some implic}, the PC-rule does not imply any of the principles in (1), 
except possibly the ED-principle.
(Also, by \ref{R+PC} if the PC-rule implies the R-principle, then it also implies the PC-principle, and
consequently also e.g. the D-principle.)
By \ref{no-implic}, the $\exists_\nabla$- and $\forall_\nabla$-principles do not imply one another.
They also do not imply any of the other principles/rules in (1) by considering models with singleton domains,
in which the $\exists_\nabla$- and $\forall_\nabla$-principles are trivially valid.

By \ref{strict implications}, the H-principle implies only the principles in (3).
By \ref{quant-no quant}, the $\forall_\Box$-principle does not imply the $\exists_\Box$-principle.
The $\forall_\Box$- and $\exists_\Box$-principles also do not imply any of the other three principles in (2) 
by considering models with singleton domains, in which the $\forall_\Box$- and $\exists_\Box$-principles 
are trivially valid.
\end{proof}

\begin{problem} \label{pc-ed}
Does the PC-rule imply the ED-principle?
\end{problem}

This question (which appeared implicitly in arXiv's version 3 of the present paper and became explicit starting from version 4)
was recently answered in the negative by A. Onoprienko \cite{On1}.

\section{Topology of Onoprienko's models}

Recently A. Onoprienko proved that QHC is complete with respect to a class of Kripke type models \cite{On2}.
In the case of the propositional fragment HC she actually found two independent proofs of completeness: 
one based on the Rasiowa--Sikorski technique with ultrafilters \cite{On0}*{\S5}, and another in a more
syntactic style \cite{On0}*{\S6}, \cite{On1}.

This section, written in December 2018, is a minor elaboration on her work.
We will describe a class of models of QHC which includes Onoprienko's models of HC as well as
a new class of topological models (not reducing to those of \S\ref{models}).

\subsection{Models with multivalued maps}

If $X$ is a topological space and $p\in X$, we will denote by $U_p$ the intersection of all 
open neighborhoods of $p$.
Let us note that if $X$ is a $T_1$ space, then $U_p=\{p\}$; but if $X$ is an Alexandrov space,
then $U_p$ is open.

We will use the language of multivalued maps.
A {\it multivalued map} $F\:X\TO Y$ is a map $F_\circ\:X\to 2^Y$ into the set of all subsets of $X$.
For each $x\in X$ we write $F(x)=F_\circ(x)$.
Every map $f\:X\to Y$ is identified with the multivalued map $F\:X\TO Y$ given by
$F(x)=\{f(x)\}$, and thus can be called a ``single-valued map''.
The {\it inverse} of a multivalued map $F\:X\TO Y$ is the multivalued map $F^{-1}\:Y\TO X$ given by
$F^{-1}(y)=\{x\in X\mid y\in F(x)\}$.
In particular, every single-valued map has a multivalued inverse.

For each $A\subset X$ we write $F(A)=\bigcup_{x\in A}F(x)$.
In particular, for each $B\subset Y$ we have $F^{-1}(B)=\{x\in X\mid F(x)\cap B\ne\emptyset\}$.
Let us note that $F^*(B):=X\but F^{-1}(Y\but B)$ consists of all $x\in X$ such that $F(x)\subset B$, and 
$\Dom(F):=F^{-1}(Y)$ consists of all $x\in X$ such that $F(x)\ne\emptyset$.

An {\it MM-structure} (MM=multivalued maps) interpreting the language of QHC is described as follows.
\begin{itemize}
\item We fix a topological space $X$, a set $S$ and a domain of discourse $\D$.
\item We further fix multivalued maps%
\footnote{The letters $i$ and $j$ may be thought of as graphic approximations of the symbols 
{\textexclamdown} and {\textquestiondown}.}
$i\:X\TO S$ and $j\:S\TO X$ such that
\begin{enumerate}
\item $\Dom(i)$ is dense in $X$;
\item $s\in i\big(j(s)\big)$ for each $s\in S$;
\item $j\big(i(p)\big)\subset U_p$ for each $p\in X$.
\end{enumerate}
\item The intuitionistic side is interpreted as in a Tarski model in open subsets of $X$.
\item The classical side is interpreted as in a Leibniz--Euler model in subsets of $S$.
\item $|\wn\Phi|=j^*(|\Phi|)$;
\item $|\oc F|=\Int\big(i^*(|F|)\big)$.
\end{itemize}

Let us note that (2) implies $\Dom(j)=S$.

\begin{theorem} \label{multi-valued}
An MM-structure is a model of QHC.
\end{theorem}

\begin{proof} 1. If $|\Phi|=X$, then $|\wn\Phi|=S$.

2. If $|F|=S$, then $|\oc F|=X$.

3. If $|F|=\emptyset$, then $|\oc F|=\emptyset$ due to (1). 
It follows that $|\neg\oc\bot|=X$.

4. We have $|\oc\wn\Phi|=\Int\Big(i^*\big(j^*(|\Phi|)\big)\Big)\supset\Int|\Phi|=|\Phi|$, the inclusion being
due to (3).
Hence $|\Phi\to\oc\wn\Phi|=X$.

5. We have $|\wn\oc F|=j^*\Big(\Int\big(i^*(|F|)\big)\Big)\subset j^*\big(i^*(|F|)\big)\subset |F|$,
the latter inclusion being due to (2).
Hence $|\wn\oc F\to F|=S$.

6. We have 
\[|\wn(\Gamma\land\Delta)|=j^*(|\Gamma\land\Delta|)=j^*(|\Gamma|\cap|\Delta|)\\
=j^*(|\Gamma|)\cap j^*(|\Delta|)=
|\wn\Gamma|\cap|\wn\Delta|=|\wn\Gamma\land\wn\Delta|.\]
Also $|\Gamma\to\Delta|=X$ implies $|\Gamma|\subset|\Delta|$, which in turn implies 
$|\wn\Gamma|\subset|\wn\Delta|$.
Hence $|\wn(\Phi\to\Psi)\land\wn\Phi|=\big|\wn\big((\Phi\to\Psi)\land\Phi\big)\big|\subset|\wn\Psi|$.
Therefore $|\wn(\Phi\to\Psi)\to(\wn\Phi\to\wn\Psi)|=|\wn(\Phi\to\Psi)\land\wn\Phi\to\wn\Psi|=S$.

7. We have 
\begin{multline*}
|\oc(P\land Q)|=\Int\big(i^*(|P\land Q|)\big)=\Int\big(i^*(|P|\cap|Q|)\big)=\Int\big(i^*(|P|)\cap i^*(|Q|)\big)\\
=\Int\big(i^*(|P|)\big)\cap\Int\big(i^*(|Q|)\big)=|\oc P|\cap|\oc Q|=|\oc(P\land Q)|.
\end{multline*}
Also $|F\to G|=S$ implies $|F|\subset|G|$, which in turn implies $|\oc F|\subset|\oc G|$.
Hence $|\oc(F\to G)\land\oc F|=\big|\oc\big((F\to G)\land F\big)\big|\subset|\oc G|$.
Therefore $|\oc(F\to G)\to(\oc F\to\oc G)|=|\oc(F\to G)\land\oc F\to\oc G|=X$.
\end{proof}

\subsection{Onoprienko's models}

A {\it special MM-structure} interpreting the language of QHC is described as follows.

\begin{itemize}
\item We fix a topological space $X$, a set $S$ and a domain of discourse $\D$.
\item We further fix multivalued maps $i\:X\TO S$ and $j\:S\TO X$ such that
\begin{enumerate}
\item[(1$'$)] $\Dom(i)=X$;
\item[(2)] $s\in i\big(j(s)\big)$ for each $s\in S$;
\item[(3)] $j\big(i(p)\big)\subset U_p$ for each $p\in X$;
\item[(4)] $i^{-1}(K)$ is closed for each $K\subset S$.
\end{enumerate}
\item The intuitionistic side is interpreted as in a Tarski model in open subsets of $X$.
\item The classical side is interpreted as in a Leibniz--Euler model in subsets of $S$.
\item $|\wn\Phi|=j^*(|\Phi|)$;
\item $|\oc F|=i^*(|F|)$.
\end{itemize}

Let us note that $i^*(|F|)=X\but i^{-1}(S\but |F|)$ is always open due to (4), so we indeed get a special case of
MM-structures.
Thus Theorem \ref{multi-valued} has the following

\begin{corollary}
Every special MM-structure is a model of QHC.
\end{corollary}

The case where $X$ is a $T_0$ Alexandrov space (in other words, a poset endowed its Alexandrov topology)
is originally due to Onoprienko \cite{On0}, \cite{On1}.
The proof of Theorem \ref{multi-valued} is not far from her argument.

\begin{theorem}[Onoprienko \cite{On0}, \cite{On1}]\label{nastia}
HC is complete with respect to the class of all special MM-models in which $X$ is a $T_0$ Alexandrov space.
\end{theorem}

\subsection{Dense image models}

Let us now consider another special case of MM-models.

A {\it dense image structure} interpreting the language of QHC is described as follows.

\begin{itemize}
\item We fix a topological space $X$, a set $S$ and a domain of discourse $\D$.
\item We further fix a map $f\:S\to X$ whose image is dense in $X$.
\item The intuitionistic side is interpreted as in a Tarski model in open subsets of $X$.
\item The classical side is interpreted as in a Leibniz--Euler model in subsets of $S$.
\item $|\wn\Phi|=f^{-1}(|\Phi|)$;
\item $|\oc F|=\Int\big(X\but f(S\but |F|)\big)$.
\end{itemize}
(Compare the interpretation of $\oc$ with generalized Tarski models of \cite{M1}*{\S\ref{int:Tarski}}.)

\begin{proposition} \label{dense is mm}
A dense image structure is an MM-structure.
\end{proposition}

\begin{proof} Let $j=f$ and $i=f^{-1}$.
Then
\begin{enumerate}
\item $\Dom(i)=f(S)$ is dense in $X$;
\item $s\in f^{-1}\big(f(s)\big)=i\big(j(s)\big)$ for each $s\in S$;
\item $j\big(i(p)\big)=f\big(f^{-1}(p)\big)=\{p\}\subset U_p$ for each $p\in X$;
\item $|\wn\Phi|=f^{-1}(|\Phi|)=j^*(|\Phi|)$;
\item $|\oc F|=\Int\big(X\but f(S\but |F|)\big)=\Int\big(i^*(|F|)\big)$.
\end{enumerate}
\end{proof}

\begin{remark} Dense image models of QHC yield the following class of models of QH4: formulas are interpreted 
by open subsets of $X$, and the modality $\nabla$ is interpreted by 
$|\nabla\Phi|=\Int\Big(X\but f\big(S\but f^{-1}(|\Phi|)\big)\Big)=\Int\Big(X\but\big(f(S)\but|\Phi|\big)\Big)$.
Let us note that this interpretation depends only on the image of $f$, and so for the purposes of interpreting
QH4 $f$ may be assumed to be an embedding.
In the case where $X$ is a $T_0$ Alexandroff space, Art\"emov and Protopopescu showed that H4 is complete 
with respect to this class of models \cite{AP}.
In fact, these models of H4 (not necessarily for $T_0$ Alexandroff spaces) are essentially well known, as 
the author learned from M. Jibladze; more precisely, well known is the situation when $\nabla\ab\to\ab$ 
is removed from the four defining laws of QH4, and accordingly the requirement that $f(S)$ be dense in $X$ is
dropped \cite{Mac}*{\S6}. 
\end{remark}

\begin{remark} Clearly, when $f=\id$, dense image models reduce to the interior-based models 
of \S\ref{interior-based}.

When $X$ is a finite (hence Alexandrov) $T_0$ space, we can also recover the regularization-based models of
\S\ref{regularization-based} from dense image models.
Namely, let $S$ be the set of all coatoms (=open singletons) of $X$, and let $f\:S\to X$ be the inclusion.
Clearly, $S$ is dense in $X$.

It is not hard to see that the boolean algebra of regular open subsets of $X$ is isomorphic to
the boolean algebra of all subsets of $S$, via the following mutually inverse maps: a regular open subset $U$
of $X$ yields the subset $U\cap S$ of $S$, and a subset $T$ of $S$ yields the regular open subset
$\Int\big(X\but(S\but T)\big)$ of $X$.
This yields a one-to-one correspondence, which preserves validity of formulas, between the dense image models 
based on $f$ and the regularization-based models based on $X$. 
\end{remark}

\begin{proposition} The following classes of models may be identified:
\begin{itemize}
\item dense image models in which $X$ is a $T_0$ Alexandrov space;
\item special MM-models in which $j$ is single-valued.
\end{itemize}
\end{proposition}

Let us note that the identification used here is different from that in Proposition \ref{dense is mm}:
$i(p)=j^{-1}(U_p)$ instead of $i=j^{-1}$.

\begin{proof}
Let us consider a special MM-model with single-valued $j$.
Then the condition $s\in i\big(j(s)\big)$ can be rewritten as $j(s)\in i^{-1}(s)$.
Since $i^{-1}(s)$ is closed, it must contain the closure of $j(s)$.
Thus $p\le j(s)$ implies $p\in i^{-1}(s)$, or equivalently $s\in i(p)$.
In other words, if $q\ge p$, then $j^{-1}(q)\subset i(p)$.
Thus $i(p)$ contains $j^{-1}(U_p)$.
But on the other hand, the condition $j\big(i(p)\big)\subset U_p$ implies that
$i(p)\subset j^{-1}(U_p)$.
Thus we have $i(p)=j^{-1}(U_p)$ for each $p\in X$.

Now $|\wn\Phi|=j^*(|\Phi|)=j^{-1}(|\Phi|)$ and $|\oc F|=i^*(|F|)=\{p\in P\mid i(p)\subset |F|\}$.
We have
\begin{multline*}
i(p)\subset |F|\iff j^{-1}(U_p)\subset |F|\iff j(S\but|F|)\subset X\but U_p\iff U_p\subset X\but j(S\but|F|)\\
\iff p\in\Int\big(X\but j(S\but|F|)\big).
\end{multline*}
Thus $|\oc F|=\Int\big(X\but j(S\but|F|)\big)$.
\end{proof}

\begin{problem} Is QHC complete with respect to the class of dense image models?
\end{problem}

The author asked A. Onoprienko this question in December 2018. 
Recently she answered it affirmatively in the propositional case \cite{On1}*{\S8}, \cite{On3}.

\begin{proposition} (a) Every MM-model in which $X$ is a $T_1$ space is equivalent to a dense image model.

(b) Every special MM-model in which $X$ is a $T_1$ space has discrete $X$.
\end{proposition}

\begin{proof}[Proof. (a)] Since $X$ is a $T_1$ space, $U_p=\{p\}$ for each $p\in X$.
Hence $j\big(i(p)\big)=\{p\}$ for each $p\in X$.
On the other hand, for each $s\in S$ we have $s\in i\big(j(s)\big)$ and in particular $s\in i(p)$ 
for some $p\in X$.
Then $j(s)\subset j\big(i(p)\big)=\{p\}$.
Thus $j$ is single-valued.
In particular, $j^*(|\Phi|)=j^{-1}(|\Phi|)$

For each $p\in X$ we have $i(p)\subset j^{-1}(p)$ due to $j\big(i(p)\big)=\{p\}$
and $j^{-1}(p)\subset i(p)$ since $j(s)=\{p\}$ implies $s\in i\big(j(s)\big)=i(p)$.
Thus $i=j^{-1}$.
In particular, $i^*(|F|)=X\but i^{-1}(S\but |F|)=X\but j(S\but |F|)$.
\end{proof}

\begin{proof}[(b)]
By the proof of (a), $j$ is single-valued and $i=j^{-1}$.
Hence $i(p)\cap i(q)=\emptyset$ whenever $p\ne q$.
Therefore for each $p\in X$ we have $i^*\big(i(p)\big)=\{p\}$.
Then $\{p\}$ is open by (4).
Thus $X$ is discrete.
\end{proof}

\subsection*{Acknowledgements}
I would like to thank A. Bauer, L. Beklemishev, M. Jibladze, A. Onoprienko, V. Shehtman and A. Shen for 
valuable discussions and useful comments.
An initial part of this work was carried out while enjoying the hospitality and
stimulating atmosphere of the Institute for Advanced Study and its
Univalent Foundations program in Spring 2013.

\subsection*{Disclaimers}

1. Some translations quoted in this paper were edited in order to improve syntactic and 
semantic fidelity.
When emphasis is present in quoted text, it is always original.

2. I oppose all wars, including those wars that are initiated by governments at the time when 
they directly or indirectly support my research. The latter type of wars include all wars 
waged by the Russian state in the last 25 years (in Chechnya, Georgia, Syria and Ukraine) 
as well as the USA-led invasions of Afghanistan and Iraq.


\begin{thebibliography}{00}

\bib{M0}{article}{
   author    = {S. A. Melikhov},
   title     = {Mathematical semantics of intuitionistic logic},
   note={\arxiv{1504.03380}},
}

\bib{M1}{article}{
   author    = {S. A. Melikhov},
   title     = {A Galois connection between classical and intuitionistic logics. I: Syntax},
   note={\arxiv{1312.2575}},
}

\bib{AP}{article}{
   author={Art{\"e}mov, S. N.},
   author={Protopopescu, T.},
   title={Intuitionistic epistemic logic},
   note={\arxiv{1406.1582v2} (not to be confused with v1, v3 or v4)},
   }

\bib{Hi1}{article}{
author={D. Hilbert},
title={Mathematische Probleme},
journal={Nachrichten von der K\"onigl. Gesellschaft der Wiss. zu G\"ottingen},
date={1900},
pages={253--297},
note={\href{http://www.digizeitschriften.de/dms/resolveppn/?PPN=GDZPPN002498863}{GDZ}.
Also published in: Archiv der Math. Phys. (3) {\bf 1} (1901), 44--63, 213--237;
English transl. {\it Mathematical problems}, Bull. Amer. Math. Soc. {\bf 8} (1902),
437--479, \href{http://www.ams.org/journals/bull/1902-08-10/S0002-9904-1902-00923-3/}{Journal}},
}

\bib{Hi2}{article}{
   author={Hilbert, D.},
   title={\"Uber das Unendliche},
   journal={Math. Ann.},
   volume={95},
   date={1926},
   pages={161--190},
   note={\href{http://resolver.sub.uni-goettingen.de/purl?GDZPPN002270641}{GDZ}.
   English transl. {\it On the infinite}, in: From Frege to G\"odel:
   A Source Book in Mathematical Logic, 1879--1931 (J. van Heijenoort, ed.), Harvard Univ. Press,
   Cambridge, MA, 1967, pp. 367--392},
}

\bib{Kol}{article}{
   author={Kolmogoroff, A.},
   title={Zur Deutung der intuitionistischen Logik},
   language={German},
   journal={Math. Z.},
   volume={35},
   date={1932},
   number={1},
   pages={58--65},
   note={\href{http://resolver.sub.uni-goettingen.de/purl?GDZPPN002373467}{GDZ};
   Russian transl. in {\it Колмогоров А. Н. Избранные труды. Математика
   и механика,} Наука, М. (1985), стр. 142--148;
   English translations: (i) {\it On the interpretation of intuitionistic logic},
   Selected Works of A. N. Kolmogorov. Vol. I,
   Mathematics and its Applications (Soviet Series), vol. 25, Kluwer,
   Dordrecht, 1991, pp. 151--158; (ii) From Brouwer to Hilbert:
   The debate on the foundations of mathematics in the 1920s (P. Mancosu, ed.),
   Oxford University Press, New York, 1998, pp. 328--334; (iii)
      \href{http://homepages.inf.ed.ac.uk/jmckinna/kolmogorov-1932.pdf}{by J. McKinna (2014)}},
}

\bib{Mac}{article}{
   author={Macnab, D. S.},
   title={Modal operators on Heyting algebras},
   journal={Algebra Universalis},
   volume={12},
   date={1981},
   number={1},
   pages={5--29},
}

\bib{On0}{thesis}{
author={A. A. Onoprienko},
title={Semantics of Kripke type for propositional logic of problems and propositions},
type={Master's Thesis (in Russian)},
organization={Moscow State University, Mech.-Math. Dept., Chair of Math. Logic and Theory of Algorithms},
date={May 2018},
}

\bib{On1}{article}{
   author={Onoprienko, A. A.},
   title={Kripke type semantics for a logic of problems and propositions},
   language={Russian},
   journal={Mat. Sb.},
   volume={211},
   date={2020},
   number={5},
   pages={98--125},
   translation={
    language={English},
    journal={Sbornik Math.},
       volume={211},
       date={2020},
       pages={709--732},
       }
}
\bib{On2}{article}{
   author={Onoprienko, A. A.},
   title={The predicate version of the joint logic of problems and propositions},
   language={Russian},
   journal={Mat. Sb.},
   volume={213},
   date={2022},
   number={7},
   pages={97--120},
   translation={
       language={English},
       journal={Sbornik Math.},
       volume={213},
       date={2022},       
       pages={to appear}
       },
}

\bib{On3}{article}{
   author={Onoprienko, A. A.},
   title={Topological models of the propositional logic of problems and propositions},
   language={Russian},
   journal={Vestnik Moskov. Univ.},
   date={2022},
   number={5},
   pages={25--30},
   translation={
          language={English},
          journal={Moscow Univ. Math. Bull.},
          date={2022},       
          pages={to appear}
          },
}

\bib{RS}{book}{
   author={Rasiowa, H.},
   author={Sikorski, R.},
   title={The Mathematics of Metamathematics},
   series={Monografie Matematyczne},
   volume={41},
   publisher={Pa\'nstwowe Wydawnictwo Naukowe},
   place={Warsaw},
   date={1963},
   translation={
   language={Russian},
   author={Е. Расёва},
   author={Р. Сикорский},
   title={Математика метаматематики},
   place={М.},
   publisher={Наука},
   date={1972},
   }
}

\bib{Sm}{article}{
author={C. Smorynski},
title={The incompleteness theorems},
book={
   title={Handbook of Mathematical Logic},
   editor={J. Barwise}
   series={Studies in Logic and the Found. of Math.},
   volume={90},
   publisher={North-Holland},
   place={Amsterdam},
   date={1977},
   },
}


\end{thebibliography}
\end{document}